\documentclass[leqno,11pt]{article}
\usepackage{amsmath,amssymb}
\usepackage{color}
\usepackage{epsfig}
\usepackage{graphicx}
\usepackage[utf8]{inputenc}
\usepackage[T1]{fontenc}
\usepackage{url}

\usepackage{dsfont}
\def\R{{\mathbb{R}}}
\def\ren{{\R^N}}

\newcounter{theorem}

\numberwithin{equation}{section}
\numberwithin{theorem}{section}
\newtheorem{theorem}{Theorem}[section]

 \newtheorem{cor}[theorem]{Corollary}
 
 \newtheorem{lem}[theorem]{Lemma}
 
 \newtheorem{prop}[theorem]{Proposition}

\newcommand{\qed}{\hfill$\square$\vspace{0.3cm}}

\begin{document}

\title{\textbf{Nonlocal operators of order near zero}}
\author{
Ernesto Correa\thanks{ernesto.correa\@@uc3m.es} \ \& Arturo de Pablo\thanks{arturop\@@math.uc3m.es (corresponding author)} \\ \\
Departamento de Matem\'{a}ticas \\ Universidad Carlos III de Madrid \\
28911 Legan\'{e}s, Spain }
\date{ }
\maketitle

\begin{abstract}
We study  Dirichlet forms defined by  nonintegrable L\'evy kernels whose singularity at the origin can be weaker than that of any fractional Laplacian. We show some properties of the associated Sobolev type spaces in a bounded domain, such as symmetrization estimates,  Hardy inequalities, compact inclusion in $L^2$ or the  inclusion in some Lorentz space. We then apply those properties to study the associated nonlocal operator $\mathfrak{L}$ and the Dirichlet and Neumann problems related to the equations $\mathfrak{L}u=f(x)$ and $\mathfrak{L}u=f(u)$ in $\Omega$.
\end{abstract}

\vskip 1cm

%
\noindent{\makebox[1in]\hrulefill}\newline

%

\section{Introduction}\label{sect-introduction}
\setcounter{equation}{0}

The aim of this paper is to study the properties of the bilinear Dirichlet form associated to a kernel $J:\mathbb{R}^N\times\mathbb{R}^N\to[0,\infty)$ and a given bounded set $\Omega\subset\mathbb{R}^N$, defined by
\begin{equation}
\label{bilinear-form}
\mathcal{E}(u,v)=\frac12\iint_{Q_\Omega}(u(x)-u(y))(v(x)-v(y))J(x,y)\,dxdy,
\end{equation}
where
\begin{equation*}
  \label{QOmega}
  Q_\Omega=(\Omega^c\times\Omega^c)^c,
\end{equation*}
and $J$ is a measurable function satisfying
\begin{equation}
\label{kernel0}
\tag{$\text{\rm H}_0$}
\left\{
\begin{array}{l}
J(x,y)\ge 0, \quad J(x,y)=J(y,x),\\[10pt]
\displaystyle
\sup_{x\in\mathbb{R}^N}\int_{\mathbb{R}^N}\min(1,| x-  y|^2)J(x,y)\,dy<\infty.
\end{array}
\right.
\end{equation}
Condition \eqref{kernel0} means that  $J$ is a L\'evy type kernel. We also assume that the kernel lies in the so-called \text{nonintegrable side}, that is,
\begin{equation}
  \label{noninteg}
\tag{$\text{\rm H}'_0$}
J(x,y)\ge \mathcal{K}(x-y)\ge0,\qquad\mathcal{K}\notin L^1(B_\varepsilon)\quad\forall\;\varepsilon>0,
\end{equation}
where $B_\varepsilon=\{z\in\mathbb{R}^N,\:|z|<\varepsilon\}$. Both these hypotheses \eqref{kernel0}--\eqref{noninteg} are assumed throughout the paper without further mention.

The power case $\mathcal{K}(z)=|z|^{-N-\alpha}$ for some $0<\alpha<2$ is well known and is related to stable processes, see also the associated operator, the fractional Laplacian,  below.  We are mostly interested in the weakly singular case which separates the fractional Laplacian side to the integrable side, i.e., $\lim_{z\to0}|z|^{N+\alpha}\mathcal{K}(z)=0$ for every $\alpha>0$.

We consider the spaces
\begin{equation}\label{H_J}
\mathcal{H}_J(\Omega)=\left\{u:\mathbb{R}^N\to\mathbb{R},\, \left.u\right|_\Omega\in L^2(\Omega),\,  \mathcal{E}(u,u)<\infty\right\}
\end{equation}
and
\begin{equation}\label{H_J0}
\mathcal{H}_{J,0}(\Omega)=\left\{u\in\mathcal{H}_J(\Omega),\, u\equiv0 \text{ in } \Omega^c\right\},
\end{equation}
with  norm
\begin{equation*}
\|u\|_{\mathcal{H}_{J}}=\left(\int_\Omega u^2+\mathcal{E}(u,u)\right)^{1/2}.
\end{equation*}

Hypothesis \eqref{kernel0} implies
$$
H_0^1(\Omega)\subset\mathcal{H}_{J,0}(\Omega)\subset\mathcal{H}_J(\Omega)\subset L^2(\Omega),
$$
if we consider the functions in $H_0^1(\Omega)$ extended by zero outside $\Omega$.
In the fractional Laplacian case $\mathcal{K}(z)=|z|^{-N-\alpha}$ for some $0<\alpha<2$ (and $N>\alpha$), then
$$
\mathcal{H}_{J}(\Omega)\subset H^{\alpha/2}(\Omega)\subset L^{\frac{2N}{N-\alpha}}(\Omega),
$$
thanks to the Hardy-Sobolev inequality, where $H^{\alpha/2}(\Omega)$ is the usual fractional Sobolev space of functions in $L^2(\Omega)$ such that $\iint_{\Omega\times\Omega}\frac{|u(x)-u(y)|^2}{|x-y|^{N+\alpha}}\,dxdy<\infty$.
As a byproduct we have $\mathcal{H}_J(\Omega)\hookrightarrow L^2(\Omega) $ compactly.

On the very other hand, in the case of integrable kernels, $\sup_{x\in\mathbb{R}^N}\int_{\mathbb{R}^N}J(x,y)\,dy= B<\infty$, (thus not satisfying~\eqref{noninteg}), we have  $\mathcal{E}(u,u)\le B\|u\|_2^2$ and therefore $\mathcal{H}_{J,0}(\Omega)\equiv L^2(\Omega)$.

A main objective of this paper is to establish the exact place were $\mathcal{H}_{J,0}(\Omega)$ and $\mathcal{H}_J(\Omega)$ lie in relation to $L^2(\Omega)$. We prove that if $\mathcal{K}$ does not oscillate too much at the origin then  the inclusion $\mathcal{H}_{J,0}(\Omega)\hookrightarrow L^2(\Omega) $ is compact. If in addition $\lim_{z\to0}|z|^N\mathcal{K}(z)=\infty$, then also $\mathcal{H}_{J}(\Omega)\hookrightarrow L^2(\Omega) $ is compact. See Theorems~\ref{th-compact} and \ref{th-compact2}.

The compacity of the inclusion $\mathcal{H}_{J,0}(\Omega)\hookrightarrow L^2(\Omega) $ can be explained by  the sharper inclusion into some Lorentz space $\mathcal{H}_{J,0}(\Omega)\hookrightarrow\mathcal{L}_{\mathcal{A},2}(\Omega)$, for some function $\mathcal{A}$ depending on $J$, see Theorem~\ref{th-Lorentz}.
In the course of the proof of this result we establish some interesting inequalities of Hardy type plus a symmetrization result. See Theorems~\ref{th-symmet},  \ref{th-hardy1} and \ref{th-hardy01}.

Associated to the bilinear form \eqref{bilinear-form} we may consider the linear operator $\mathfrak{L}: \mathcal{H}_J(\Omega)\to \mathcal{D}'(\Omega)$  that satisfies
$\langle \mathfrak{L}u,\zeta\rangle =\mathcal{E}(u,\zeta)$ for any $u,\zeta\in C^\infty_{\rm 0}(\mathbb{R}^N)$, where $\langle\cdot,\cdot\rangle$ denotes the dual product in the space of distributions. In a very general setting the operator $\mathfrak{L}$ is defined for all $u\in C^2(\mathbb{R}^N)$ satisfying some growth condition at infinity by
\begin{equation}
\label{operatorL}\mathfrak{L}u(x)=\textrm{P.V.}
\int_{\mathbb{R}^{N}}(u(x)-u(y))J(x,y)\,dy.
\end{equation}
When $J(x,y)=|x-y|^{-N-\sigma}$, the operator is a multiple of the well known Fractional Laplacian
\begin{equation}
\label{Fractlap}(-\Delta)^{\alpha/2}u(x)=C_{N,\alpha}\textrm{P.V.}
\int_{\mathbb{R}^{N}}\frac{u(x)-u(y)}{|x-y|^{N+\alpha}}\,dy,
\end{equation}
where $C_{N,\alpha}$ is some normalizing constant. This is a pseudo-differential operator of order $\alpha$, and behaves like $\alpha$ derivatives.
On the other hand, if $J(x,y)=\mathcal{K}(x-y)$, and $\mathcal{K}\in L^1(\mathbb{R}^N)$, with for instance $\int_{\mathbb{R}^N}\mathcal{K}=1$, then $\mathfrak{L}$ is given by
\begin{equation*}
\label{JdeJ}
\mathfrak{L}u=u-\mathcal{K}*u.
\end{equation*}
This implies that $\mathfrak{L}$ is a zero order operator, since $u$ and $\mathfrak{L}u$ possesses  the same regularity.

These two types of operators give rise to two lines of research, often disconnected. The threshold between those types of operators  is what motivates this work. We note that letting $\alpha\to0^+$ in~\eqref{Fractlap} we obtain the identity operator, precisely because of the normalizing constant, since $C_{N,\alpha}\sim\alpha\to0^+$. What we want to study here, in some very informal way, is  the limit $\alpha\sim0^+$ in the singularity of the kernel near the origin, but without the normalizing constant. We show some properties of the operator $\mathfrak{L}$ in the limit of singularity, and study the effect of applying $\mathfrak{L}$ to continuous functions, showing that the operator $\mathfrak{L}$ can be considered of order near zero but positive, see Corollary~\ref{sharp-regularity}. We also study when $\mathfrak{L}\mathds{1}_\Omega$ is integrable in $\Omega$, leading to the concept of $J$-perimeter, Theorem~\ref{perimeter}.

We also study the eigenvalues of  $\mathfrak{L}$ in $\Omega$, that is the equation $\mathfrak{L}\varphi=\lambda\varphi$ for $\varphi\in\mathcal{H}_{J,0}(\Omega)$. Existence follows from the compact inclusion $\mathcal{H}_{J,0}(\Omega)\subset L^2(\Omega)$, and then we characterize the space $\mathcal{H}_{J,0}(\Omega)$ by means of the eigenvalues, see Theorem~\ref{prop-eigen} and Proposition~\ref{eigen-op}. We also estimate the first eigenvalue in terms of the size of $\Omega$, Theorem~\ref{berezin}.

We then pass to study the problem
$$
\begin{cases}
\mathfrak{L}u=f,& \mbox{in }\Omega,\\ u=0,& \mbox{in }\Omega^c.
\end{cases}
$$
Existence and uniqueness is easily established for $f\in H^*(\Omega)$, the dual of $\mathcal{H}_{J,0}(\Omega)$.
We are  interested in the regularizing effects, and prove that $u$ has slightly better integrability than $f$ in the sense that if $f\in L^p(\Omega)$, $p\ge2$,  then $u\in \mathcal{L}_{\mathcal{A},p}(\Omega)$, a Lorentz type space. See Theorem~\ref{th-Lp-smoothing}.

The next task is to consider nonlinear problems of the form
$$
\begin{cases}
\mathfrak{L}u=f(u),& \mbox{in }\Omega,\\ u=0,& \mbox{in }\Omega^c.
\end{cases}
$$
We prove existence and uniqueness when $f$ is sublinear, see Theorem~\ref{th-sublinear}. As for superlinear reactions we show that  there is no solution when $f(u)=u^p$, $p>\frac{N+\sigma}{N-\sigma}$ and $\Omega$ is star-shaped, where $\sigma$ depends on the kernel $J$, Corollary~\ref{cor-p^*}. The proof is based on a Pohozaev inequality as obtained by~\cite{RosOtonSerra} in the fractional Laplacian case.

We finally  study  a  Neumann problem associated to the operator $\mathfrak{L}$ as it is done in \cite{diPierro-RosOton-Valdinoci} for the fractional Laplacian,
$$
\begin{cases}
\mathfrak{L}u=f,& \mbox{in }\Omega,\\ \mathcal{N}u=0,& \mbox{in }\Omega^c,
\end{cases}
$$
where $\mathcal{N}$ is some operator generalizing the normal derivative, and show some properties of the solution in our situation. See Theorem~\ref{th-Neumann-exist}.

\subsection{Organization of the paper}

We devote a preliminary Section~\ref{sect-preliminaries} to settle the precise hypotheses that we consider throughout the paper and to prove the compact inclusions of our Sobolev type spaces in $L^2$; we also present a symmetrization result. In Section~\ref{sect-functional} we prove two Hardy inequalities. The inclusion of the Sobolev spaces into some space of Lorentz type is shown in Section~\ref{sect-lorentz}. Section~\ref{sect-eigenvalues} is devoted to study the associated nonlocal operator and its properties, such as the existence of eigenvalues or its action on different functions. In Sections~\ref{sect-problems}, \ref{sect-nonlinear} and \ref{sect-neumann} we study three problems associated to that operator, two linear problems, with Dirichlet or Neumann exterior condition, and a nonlinear problem with different reactions; we show existence and uniqueness for sublinear reactions and nonexistence when the reaction is supercritical.

\section{Preliminaries}\label{sect-preliminaries}
\setcounter{equation}{0}

In order to describe the properties of the fractional type spaces $\mathcal{H}_J$ defined in \eqref{H_J} and \eqref{H_J0}, in terms of the kernel, we write condition~\eqref{noninteg} in the form
\begin{equation}\label{kernel2}\tag{$\text{\rm H}_1$}
\left\{
\begin{array}{l}
  \mathcal{K}(z)= |z|^{-N}\ell(|z|)\qquad\text{for } 0<|z|<\rho, \\[10pt]
\displaystyle M(r):=\int_r^\rho\frac{\ell(s)}s\,ds\to\infty\quad \text{as } r\to0^+,
\end{array}
\right.
\end{equation}
for some function $\ell:(0,\rho)\to(0,\infty)$, $\,\rho>0$ satisfying $0<c_1(\varepsilon)\le \ell(s)\le c_2(\varepsilon)<\infty$ for every $0<\varepsilon<\rho$. We establish some inequalities in terms of the function $M$.
We first observe that the bilinear form \eqref{bilinear-form}, when applied to functions vanishing outside $\Omega$, coincides with the global bilinear form
\begin{equation*}
\mathcal{E}_1(u,u)=\frac12\int_{\mathbb{R}^N}%
\int_{\mathbb{R}^N}|u(x)-u(y)|^2J(x,y)\,dxdy,
\end{equation*}
which is adequate also to study problems defined in the whole space. But for problems defined in a bounded domain, with a nontrivial condition in the complement of the domain, the associated bilinear form is $\mathcal{E}$ and not $\mathcal{E}_1$. See Subsection~\ref{sub-nonhom}.

On the other hand, $\mathcal{E}$  is also different from the one  related with the so-called censored processes,
\begin{equation*}
\mathcal{E}_c(u,u)=\frac12\int_{\Omega}%
\int_{\Omega}|u(x)-u(y)|^2J(x,y)\,dxdy,
\end{equation*}
see for instance \cite{BogdanBurdzyChen}. Actually,
\begin{equation}\label{otraE}
\mathcal{E}(u,u)=\mathcal{E}_c(u,u)+\int_{\Omega}|u(x)|^2\Lambda(x)\,dx,
\end{equation}
where
\begin{equation}\label{Lambda}
\Lambda(x)=\int_{\Omega^c}J(x,y)\,dy.
\end{equation}
Clearly the second integral in \eqref{otraE} is strictly positive. If for instance  $J(x,y)\ge c>0$ for every $|x-y|\le R$, and $R>\delta=\sup_{x\in\Omega}dist(x,\partial\Omega)$  then $\Lambda(x)\ge c|\{\delta<|x-y|<R\}|=A>0$ for every $x\in\Omega$.
See Theorem~\ref{th-hardy01} below for a more precise estimate of this function. This gives the Poincar\'e inequality
\begin{equation}\label{poincare}
\mathcal{E}(u,u)\ge A\|u\|_2^2,
\end{equation}
or which is the same, the property
\begin{equation}\label{eq-L2}
\mathcal{H}_{J,0}(\Omega)\subset L^2(\Omega).
\end{equation}
Even more, the Poincar\'e inequality holds also for the bilinear form $\mathcal{E}_c$, though the proof is not so direct.
A very much weaker condition to have a Poincar\'e inequality is obtained in \cite{Felsinger-Kassmann-Voigt}, where they prove that it is enough to have $|\{K(z)>0\}|>0$.

To show that the inclusion \eqref{eq-L2} is proper, and also compact, for non integrable kernels, we  impose that the function $\ell$ in \eqref{kernel2}  vary slowly at the origin, that is
\begin{equation}\label{slowly}\tag{$\text{\rm H}_2$}
  \lim_{s\to0}\frac{\ell(\lambda s)}{\ell(s)}=1\qquad\text{for every } \lambda>0.
\end{equation}
See the monograph \cite{BGT} for the properties of slowly varying functions.
Examples of slowly varying functions $\ell$ that also satisfy \eqref{kernel2} are $\ell(s)=1$, $\ell(s)=\log^\beta(2\rho/s)$, $\beta\ge-1$, or $\ell(s)=(\log(2\rho/s)\log\log(2\rho/s))^{-1}$. We do not consider highly oscillating functions like $\ell(s)=1-\sin(1/|z|)$. We may also assume, without loss of generality, that $\ell$ is differentiable near the origin and
\begin{equation}\label{slowly2}
  \lim_{s\to0}\frac{s\ell'(s)}{\ell(s)}=0,
\end{equation}
since every slowly varying function can be controlled from below by another slowly varying function with that property, see the representation formula \cite[Theorem 1.3.1]{BGT}.

\subsection{Compact embeddings in $L^2$}\label{sect-inclusion}

\begin{theorem}
  \label{th-compact} In the above hypotheses \eqref{kernel2}--\eqref{slowly} the embedding $\mathcal{H}_{J,0}(\Omega)\hookrightarrow L^2(\Omega)$ is compact.
\end{theorem}
In order to prove the same property for the bigger space $\mathcal{H}_J(\Omega)$ we must add an extra hypothesis.

\begin{theorem}\label{th-compact2}
Assume  hypothesis \eqref{kernel2} with $\ell(0^+)=\infty$. Then the embedding $\mathcal{H}_J(\Omega)\hookrightarrow L^{2}(\Omega)$ is compact.
\end{theorem}

For the first result we consider the bilinear form $\mathcal{E}_{\mathcal{K}}$ associated to the convolution kernel $\mathcal{K}$ of~\eqref{kernel2}, that is
$$
\mathcal{E}_{\mathcal{K}}(u,u)=\frac12\int_{\mathbb{R}^{N}}%
\int_{\mathbb{R}^{N}}(u(x)-u(y))^2\mathcal{K}(x-y)\,dxdy.
$$
Recall that $\mathcal{E}_1(u,u)\ge\mathcal{E}_{\mathcal{K}}(u,u)$. On the other hand
$$
\mathcal{E}_{\mathcal{K}}(u,u)=\int_{\mathbb{R}^N}m(\xi)|u(\xi)|^2\,d\xi,
$$
where the multiplier $m$ is given by
$$
m(\xi)=\int_{\mathbb{R}^N} \big (1-\cos(z\cdot \xi)\big)\mathcal{K}(z)\,d z.
$$
This multiplier has been estimated in \cite{Kassmann-Schwab} using the function $M(r)$ and hypotheses \eqref{slowly},
\begin{equation*}
m(\xi)\ge cM(|\xi|^{-1}), \qquad\text{for every } |\xi|>1.
\end{equation*}

We are now in a position to prove Theorem~\ref{th-compact}, which is a direct application of the following characterization of Pego, see \cite{Pego}.

\begin{theorem}\label{th-pego}
A bounded subset $\Sigma$ of $L^2(\mathbb{R}^N)$ is conditionally compact if and only if
\begin{equation}\label{compactPego}
\lim_{R\to\infty}\sup_{f\in \Sigma}\int_{|x|>R}|f(x)|^2\,dx=\lim_{R\to\infty}\sup_{f\in \Sigma}\int_{|\xi|>R}|\widehat{f}(\xi)|^2\,d\xi=0.
\end{equation}
\end{theorem}

\noindent{\it Proof of Theorem~\ref{th-compact}.}
For a constant $C>0$ let
$$
\Sigma=\{f\in \mathcal{H}_{J,0}(\Omega)\,:\, \|f\|_{\mathcal{H}_J}\le C\}\subset L^2(\mathbb{R}^N).
$$
We first have, since $\Omega$ is bounded, that if $R$ is large enough
$$
\int_{|x|>R}|f(x)|^2\,dx=0, \qquad\text{ for every } f\in\Sigma.
$$
On the other hand, from the previous calculations we have,
$$
C^2\ge \mathcal{E}(f,f)=\mathcal{E}_1(f,f)\ge\mathcal{E}_{\mathcal{K}}(f,f)\ge c\int_{|\xi|>R}M(|\xi|^{-1})|\widehat{f}(\xi)|^2\,d\xi.
$$
Thus
$$
\int_{|\xi|>R}|\widehat{f}(\xi)|^2\,d\xi\le \frac c{M(1/R)}, \qquad\text{ for every } f\in\Sigma,
$$
since $M$ is nonincreasing. We conclude with the fact that $M(0^+)=\infty$ that \eqref{compactPego} holds.
\qed

We now consider the sharper compact inclusion $\mathcal{H}_J(\Omega)\hookrightarrow L^{2}(\Omega)$. It requires the extra hypothesis $\ell(0^+)=\infty$.

\noindent{\it Proof of Theorem \ref{th-compact2}.}
The proof follows the one of the classical Riesz-Fr\'echet-Kolmogorov Theorem, as adapted in \cite{PalatucciSavinValdinoci} for the fractional Laplacian.

Let $\mathcal{F}\subset{\mathcal{H}}_J(\Omega)$ be a bounded set. We show that $\mathcal{F}$ is totally bounded in $L^2(\Omega),$ i.e., for any $\epsilon \in (0,1)$ there exist $\beta_1,...,\beta_M \in L^2(B_1)$ such that for any $u\in \mathcal{F}$ there exists $j\in\{1,...,M\}$ such that
\begin{equation*}\label{H131}
\|u-\beta_j\|_{L^2(\Omega)}\leq \epsilon.
\end{equation*}
We take a collection of disjoints cubes $Q_1,...Q_R$ of side $\rho<1$ such that $\Omega=\bigcup_{j=1}^{R} Q_j$.
For any $x\in \Omega$ we define $j(x)$ as the unique integer in $\{1,...,R\}$ for which $x\in Q_{j(x)}$.
Also, for any $u \in \mathcal{F},$ let
$$
P(u)(x):=\frac{1}{|Q_{j(x)}|}\int_{Q_{j(x)}}u(y)\,dy.
$$
Notice that
$$
P(u+v)=P(u)+P(v) \;\mbox{for any}\; u,v\in \mathcal{F},
$$
and that $P(u)$ is constant, say equal to $q_j(u)$, in any $Q_j$, for $j\in\{1,...,R\}$.
Therefore, we can define
$$
S(u):=\rho^{N/2}\left(q_1(u),...,q_R(u)\right)\in \ren.
$$
We observe that $S(u+v)=S(u)+S(v)$. Moreover,
\begin{equation}\label{H133}
\begin{array}{rl}
\|P(u)\|_{L^2(\Omega)}^{2}&\displaystyle=\sum_{j=1}^R\int_{Q_j}|P(u)|^2\,dx\\
&\displaystyle\leq \rho^N\sum_{j=1}^R|q_j(u)|^2=|S(u)|^2\leq\frac{|S(u)|^2}{\rho^N}.
\end{array}
\end{equation}
And, by H\"older inequality,
\begin{equation*}
\begin{array}{rl}
|S(u)|^2&\displaystyle=\sum_{j=1}^{R}\rho^N|q_j(u)|^2=\frac{1}{\rho^N}\sum_{j=1}^{R}\left|\int_{Q_j}u(y)\,dy\right|^2\\
&\displaystyle\leq \sum_{j=1}^{R}\int_{Q_j}|u(y)|^2\,dy=\int_{\Omega}|u(y)|^2\,dy=\|u\|_{L^2(\Omega)}^2.
\end{array}
\end{equation*}
In particular,
$$
\sup_{u\in\mathcal{F}}|S(u)|^2\leq C,
$$
that is, the set $S(\mathcal{F})$ is bounded in $\ren$ and so, since it is finite dimensional, it is totally bounded. Therefore, there exist $b_1,...,b_K\in\ren$ such that
\begin{equation}\label{H134}
S(\mathcal{F})\subset\bigcup_{i=1}^{K}B_{\eta}(b_i).
\end{equation}
For any $i\in\{1,...,K\}$, we write the coordinates of $b_i$ as $b_i=(b_{i,1},...,b_{i,N})\in \ren$. For any $x\in \Omega$, we set
$$
\beta_i(x):=\rho^{-N/2}b_{i,j(x)},
$$
where $j(x)$ is as above.
Notice that $\beta_i$ is constant on $Q_j$, i.e. if $x\in Q_j$ then
\begin{equation}\label{H135}
P(\beta_j)(x)=\rho^{-N/2}b_{i,j}=\beta_{i}(x)
\end{equation}
and so $q_j(\beta_i)=\rho^{-N/2}b_{i,j}$; thus $S(\beta_i)=b_i$.
Furthermore, for any $u\in \mathcal{F}$, by H\"older inequality,
\begin{equation*}
\begin{array}{l}
\|u-P(u)\|_{L^2(\Omega)}^{2}\displaystyle=\sum_{j=1}^{R}\int_{Q_j}\left|u(x)-P(u)(x)\right|^2\,dx=\sum_{j=1}^{R}\int_{Q_j}\left|u(x)-\frac{1}{|Q_j|}\int_{Q_j}u(y)\,dy\right|^2\,dx\\
\qquad\displaystyle=\frac{1}{\rho^{2N}}\sum_{j=1}^{R}\int_{Q_j}\left|\int_{Q_j}\left(u(x)-u(y)\right)\,dy\right|^2\,dx
\le\frac{1}{\rho^{N}}\sum_{j=1}^{R}\int_{Q_j}\int_{Q_j}\left|u(x)-u(y)\right|^2\,dy\,dx\\
\qquad\displaystyle\leq\frac{1}{\ell(\rho)}\sum_{j=1}^{R}
\int_{Q_j}\int_{Q_j}\left|u(x)-u(y)\right|^2J(x-y)\,dy\,dx\leq\frac{2{\mathcal{E}}(u,u)}{\ell(\rho)}.
\end{array}
\end{equation*}
Consequently, for any $j\in\{1,...,K\}$, recalling \eqref{H133} and \eqref{H135}
\begin{equation}\label{H137}
\begin{array}{rl}
\|u-\beta_j\|_{L^2(\Omega)}&\displaystyle\leq\|u-P(u)\|_{L^2(\Omega)}+
\|P(\beta_j)-\beta_j\|_{L^2(\Omega)}+\|P(u-\beta_j)\|_{L^2(\Omega)}\\ [3mm]
&\displaystyle\leq \frac{2{\mathcal{E}}(u,u)}{\ell(\rho)}+
\frac{\left|S(u)-S(\beta)\right|}{\rho^{N/2}}.
\end{array}
\end{equation}
Now, given any $u\in \mathcal{F}$, we recall \eqref{H134} and  we take $j\in \{1,...,K\}$ such that $S(u)\in B_{\eta}(b_j)$. Then, \eqref{H135} and \eqref{H137} give that
$$
\|u-\beta_j\|_{L^2(\Omega)}\leq \frac{2{\mathcal{E}}(u,u)}{\ell(\rho)}+\frac{\left|S(u)-S(\beta)\right|}{\rho^{N/2}}\leq \frac{2{\mathcal{E}}(u,u)}{\ell(\rho)}+\frac{\eta}{\rho^{N/2}}<\epsilon.
$$
\qed

We also recall the work \cite{ChenHajaiej}, where the authors consider symmetric kernels  also in the limit of integrability. But their condition
$$
\sup\left\{s\ge0\,:\,\lim_{r\to0^+}r^{s}\int_{|x-y|>r}J(x-y)\,dy=\infty\right\}>0,
$$
is not fulfilled in general by the kernels studied in this paper. In particular it requires $\limsup_{|x|\to0^+}|x|^{N+\varepsilon}J(x)>0$ for some $\varepsilon>0$.

\subsection{Symmetrization}

An easy property of bilinear forms of the type \eqref{bilinear-form} is that they decrease when taking absolute values
\begin{equation}\label{absolute}
\mathcal{E}_1(|u|,|u|)\le \mathcal{E}_1(u,u).
\end{equation}
This follows from the inequality $\left||a|-|b|\right|\le |a-b|$. It is also a consequence of the following more general inequality, called Stroock-Varopoulos  inequality \cite{Varopoulos,BrandledePablo}, which will be used later on.
 \begin{prop}\label{Stroock_Varo}
 	Let $u\in \mathcal{H}_J(\mathbb{R}^N)$ such that $F(u),\, G(u),\,\Phi(u)\in \mathcal{H}_J(\mathbb{R}^N)$, and assume $ \left(\Phi'\right)^2\le F'G'$. Then
 	\begin{equation*}
\mathcal{E}_1(\Phi(u),\Phi(u))\le\mathcal{E}_1(F(u),G(u)).
 	\end{equation*}	
 \end{prop}
Clearly the same is true with $\mathcal{E}_1$ replaced by $\mathcal{E}$ provided $u\in \mathcal{H}_{J,0}(\mathbb{R}^N)$. It is not clear what happens in the last case for general $u\in \mathcal{H}_J(\mathbb{R}^N)$.

We prove in this section that the energy $\mathcal{E}(u,u)$ also decreases when we replace $u$ by its  symmetric rearrangement, provided the kernel is radially symmetric and decreasing. This property is well known for the norm in $H_0^{\alpha/2}(\Omega)$, $0<\alpha\le2$.

For a measurable function we consider its distribution function
$$
\mu(t)=|\{x\in\mathbb{R}^N:|u(x)|>t\}|.
$$
We then define the decreasing rearrangement $u^*$ of $u$ to be the radially deacreasing function with the same distribution function as $u$, that is, $u^*:B_R\to\mathbb{R}^+$, $R=\left(\frac{|\Omega|}{\omega_N}\right)^{1/N}$, satisfies
$$
u^*(x)=\inf\{\lambda>0\,:\,\mu(\lambda)<\omega_N|x|^N\}.
$$
Here $\omega_N$ is the measure of the unit ball in $\mathbb{R}^N$. We also have the layer cake representation
\begin{equation}\label{layercake}
u^*(x)=
\int_0^\infty \mathds{1}_{\{u(x)> t\}}\,dt,
\end{equation}
and
$$
\int_{B_R}u^*(x)\,dx=\int_\Omega|u(x)|\,dx=\int_0^\infty\mu(t)\,dt.
$$

All the $L^p$ norms are conserved as well under symmetrization, or even the integral of any convex, nonnegative symmetric function of $u$. On the other hand, for Sobolev norms symmetrization decreases, that is
$$
\int_{B_R}|\nabla u^*(x)|^p\,dx\le\int_\Omega|\nabla u(x)|^p\,dx,
$$
or even for fractional Sobolev norms
$$
\int_{\mathbb{R}^N}|(-\Delta)^{\alpha/2} u^*(x)|^2\,dx\le\int_{\mathbb{R}^N}|(-\Delta)^{\alpha/2} u(x)|^2\,dx.
$$
See for instance \cite{Almgren-Lieb}. We prove here that symmetrization also decreases the norm in the space $\mathcal{H}_J(\Omega)$.

\begin{theorem}
  \label{th-symmet} Assume $J(x,y)=\mathcal{K}(|x-y|)$ where $\mathcal{K}$ is a nonincreasing  function that satisfies   \eqref{kernel2}. If $u\in\mathcal{H}_{J,0}(\Omega)$ and $u^*$ is its decreasing rearrangement, then
  \begin{equation*}\label{symm}
    \mathcal{E}(u,u)\ge \mathcal{E}(u^*,u^*).
  \end{equation*}
\end{theorem}

We use the following result, which is the key point to prove the inequality in the fractional framework, see again \cite{Almgren-Lieb}.

\begin{theorem}\label{th-A-L}
Let $\Phi\in L^1(\mathbb{R}^N)$
be a positive symmetric decreasing function. Then for all non-negative measurable $u$ one has
$$
\int_{\mathbb{R}^{N}}%
\int_{\mathbb{R}^{N}}(u(x)-u(y))^2\Phi(x-y)\,dxdy\ge
\int_{\mathbb{R}^{N}}%
\int_{\mathbb{R}^{N}}(u^*(x)-u^*(y))^2\Phi(x-y)\,dxdy.
$$
\end{theorem}

\noindent{\it Proof of Theorem \ref{th-symmet}.}
We write, as in \cite{Almgren-Lieb},
$$
\mathcal{E}(u,u)
=\frac12\int_{\mathbb{R}^{N}}%
\int_{\mathbb{R}^{N}}|u(x)-u(y)|^2\left(\int_0^\infty e^{\frac{-t}{\mathcal{K}(|x-y|)}}\,dt\right)\,dxdy
=\int_0^\infty G(u,t)\,dt
$$
where
$$
G(u,t)=\frac12\int_{\mathbb{R}^{N}}%
\int_{\mathbb{R}^{N}}|u(x)-u(y)|^2\Phi(|x-y|)\,dxdy,\qquad\Phi(z)=e^{\frac{-t}{\mathcal{K}(z)}}.
$$
By adding $\varepsilon$ times a positive integrable function to our kernel we get that $\Phi\in L^1(\mathbb{R}^N)$ satisfies the hypotheses of Theorem~\ref{th-A-L}. We then conclude, by letting $\varepsilon\to0$,
$$
G(u,t)\ge G(u^*,t)\,\qquad\text{for every } t>0,
$$
and thus
$$
\mathcal{E}(u,u)\ge \mathcal{E}(u^*,u^*).
$$
\qed

\section{Hardy inequalities}\label{sect-functional}
\setcounter{equation}{0}

We now establish some interesting inequalities for the bilinear form $\mathcal{E}$ which will be used afterwards in order to sharpen the inclusion~\eqref{eq-L2}.
We prove two Hardy inequalities of the form
$$
    \mathcal{E}(u,u)\ge c\int_{\Omega} |u(x)|^2\psi(x)\,dx,
$$
where $\psi$ is a weight that can be singular at the origin or at the boundary of $\Omega$. See for instance \cite{Opic-Kufner} for the classical Hardy  inequalities in the local case, and \cite{Dyda,Frank-Seiringer,Herbst} for the fractional Laplacian case.

We begin with a Hardy inequality that contains  a weight singular at the origin. It is to be compared with the classical Hardy inequality,
\begin{equation}\label{hardy-lap}
\int_{\Omega}|\nabla u(x)|^2\,dx\ge d_N\int_{\Omega} \frac{|u(x)|^2}{|x|^2}\,dx,
\end{equation}
and
the fractional Hardy inequality, corresponding to $J(x,y)=|x-y|^{-N-\alpha}$, $0<\alpha<2$, which is
\begin{equation}\label{hardy-lapfrac}
\int_{\mathbb{R}^{N}}%
\int_{\mathbb{R}^{N}}\frac{|u(x)-u(y)|^2}{|x-y|^{N+\alpha}}\,dxdy\ge d_{N,\alpha}\int_{\Omega} \frac{|u(x)|^2}{|x|^\alpha}\,dx,
\end{equation}
for some explicit constants $d_N$ and $d_{N,\alpha}$.
In our situation of weakly singular kernels the weight depends on the function $\ell$, and if  for instance $\ell(0)=c>0$ it is logarithmic. In the proof we use an estimate of the action of the nonlocal operator $\mathfrak{L}$ in~\eqref{operatorL} over some power. We postpone this estimate to Section~\ref{sect-eigenvalues} where we study the properties of that operator.

\begin{theorem}
  \label{th-hardy1} Assuming hypotheses  \eqref{kernel2}--\eqref{slowly},   for every $u\in\mathcal{H}_{J,0}(\Omega)$ it holds
  \begin{equation}\label{hardy1}
    \mathcal{E}(u,u)\ge c\int_{\Omega} |u(x)|^2 M(\rho|x|/R)\,dx,
  \end{equation}
  $R=\sup_{x\in\Omega}|x|$. If moreover $\ell(0)>0$ then
\begin{equation*}\label{hardy-log}
    \mathcal{E}(u,u)\ge c\int_{\Omega} |u(x)|^2\left|\log(\rho|x|/R)\right|\,dx,
\end{equation*}
\end{theorem}
\noindent{\it Proof.}
We first observe that if $w$ is a nontrivial nonnegative function, then
$$
\begin{array}{rl}
\displaystyle\mathcal{E}(u,u)&\displaystyle=\mathcal{E}\left(\frac{u^2}w,w\right) +
\frac12\int_{\mathbb{R}^{N}}%
\int_{\mathbb{R}^{N}}w(x)w(y)\left(\frac{u(x)}{w(x)}-\frac{u(y)}{w(y)}\right)^2J(x,y)\,dxdy,\\ [4mm]
&\displaystyle\ge\int_\Omega\frac{u^2}{w}\,\mathfrak{L}w
\end{array}
$$
since the second integral on the right is nonnegative. Now put $w(x)=|x|^{-N/2}$ and use Lemma~\ref{lema.x-N} to get
$$
\mathcal{E}(u,u)\ge c\int_{|x|<\varepsilon} |u(x)|^2M(|x|)\,dx.
$$
With Poincar\'e inequality this gives \eqref{hardy1}, since $M$ is bounded outside the origin, and $M$ is rescaled in order to be defined at all points of $\Omega$.
\qed

{\bf Remarks.}
$i)$ The first equality in the proof is a sort of Picone identity \cite{Picone}, to be compared with the classical one
$$
|\nabla u|^2=\nabla\left(\frac{u^2}w\right)\cdot\nabla w+u^2\left|\frac{\nabla u}u-\frac{\nabla w}w\right|^2,
$$
provided $u,\,w\in C^1(\Omega)$, $w\ge0$, $w\not\equiv0$.

$ii)$ Observe that  what we indeed have is a Hardy inequality with remainder
$$
\mathcal{E}(u,u)\ge c\int_{\Omega} |u(x)|^2M(|x|)\,dx+\frac12\int_{\mathbb{R}^{N}}%
\int_{\mathbb{R}^{N}}\frac{\left(|x|^{N/2}u(x)-|y|^{N/2}u(y)\right)^2}{|x|^{N/2}|y|^{N/2}}J(x,y)\,dxdy.
$$

$iii)$ In general we obtain inequality \eqref{hardy1} for a given weight function $\psi$ if there exists a function $w$ such that $\mathfrak{L}w\ge\psi w$.

$iv)$ One of the main features of the classical Hardy inequalities~\eqref{hardy-lap} and~\eqref{hardy-lapfrac} is that the optimal constants can be obtained in a precise way. And this depends on the fact that we can obtain an explicit function $w$ for which $\mathfrak{L}w=\psi w$, see the comment after Lemma~\ref{lema.x-N}. This is not possible in our situation.

We also obtain a Hardy inequality in a ball with a weight singular at the boundary, in the spirit of~\cite{Dyda}.

\begin{theorem}\label{th-hardy01}
Assuming hypotheses  \eqref{kernel2}--\eqref{slowly},  for every $u\in\mathcal{H}_{J,0}(B_1)$ it holds
  \begin{equation*}\label{hardy01}
    \mathcal{E}(u,u)\ge  c\int_{B_1} |u(x)|^2M(1-|x|)\,dx.
  \end{equation*}
\end{theorem}
\noindent{\it Proof.}
By \eqref{otraE} we only have to estimate the function $\Lambda(x)=\int_{\{|x-y|>1,\,|y|<1\}}|y|^{-N}\ell(|y|)\,dy$.
For each $|x|<1$ given, let $R=R_x$ be the rotation that carries $x$ to the negative first  axis, that is, $w=Rx=-|x|e_1$, where $e_1=(1,0,\dots,0)$, and perform the change of variables  $z=Ry$. Since
$$
z_1+|x|=z_1+|w|>1\implies z_1-w_1>1\implies |z-w|>1,
$$
we get, if $N\ge2$,
$$
\begin{array}{rl}
\Lambda(x)&\displaystyle\ge c\int_{1-|x|}^1\int_0^1(t^2+\rho^2)^{-N/2}\rho^{N-2}\ell\left(\sqrt{t^2+\rho^2}\right)\,d\rho dt \\ [4mm]
&\displaystyle\ge c\int_{1-|x|}^1\int_0^1 t^{-1}w^{N-2}\ell\left(t\sqrt{1+w^2}\right)\,dw dt.
\end{array}
$$
Property \eqref{slowly} implies
$$
\ell\left(t\sqrt{1+w^2}\right)\ge c\ell(t)\qquad\text{for every } 0<w<1,
$$
so that
$$
\Lambda(x)\ge c\int_{1-|x|}^1 \int_0^1\frac{\ell(t)}t \,w^{N-2}\,dw dt=cM(1-|x|).
$$
If $N=1$ we get the same estimate directly.
\qed

\section{Lorentz spaces}\label{sect-lorentz}

For a given function $\mathcal{A}:\mathbb{R}^+\to\mathbb{R}^+$, increasing with $\mathcal{A}(0)=0$, and a constant $p\ge1$, we define the Lorentz-type space
$$
 \mathcal{L}_{\mathcal{A},p}(\mathbb{R}^N)=\left\{u \text{ measurable, }  \int_0^\infty\mathcal{A}(\mu(t))t^{p-1}\,dt<\infty\right\},
$$
with seminorm
$$
\|u\|_{\mathcal{A},p}=\left(p\int_0^\infty\mathcal{A}(\mu(t))t^{p-1}\,dt\right)^{1/p},
$$
where $\mu(t)$
is the distribution function of $u$.
We may replace $\mathbb{R}^N$ by $\Omega$ in the  definitions if we restrict ourselves to functions $u$  that vanish outside $\Omega$. In that case the above is a norm. These spaces generalize the Lebesgue spaces $L^p$ and the standard Lorentz spaces $L_{q,p}$. In fact
$$\mathcal{L}_{\mathcal{A},p}(\Omega)= \left\{
\begin{array}{ll}
L_{q,p}(\Omega)&\quad\text{if } \mathcal{A}(s)=s^{p/q}\\ [2mm]
L^{p}(\Omega)&\quad\text{if } \mathcal{A}(s)=s.
\end{array}\right.
$$

\begin{theorem}
  \label{th-Lorentz} Assume hypotheses  \eqref{kernel2}--\eqref{slowly} with $\mathcal{K}$ a radially symmetric nonincreasing  function. Then there exists a function $\mathcal{A}:[0,|\Omega|]\to\mathbb{R}^+$ such that  $\mathcal{H}_{J,0}(\Omega)\subset \mathcal{L}_{\mathcal{A},2}(\Omega)$.
\end{theorem}

The function $\mathcal{A}$ depends on $\mathcal{K}$ through formula~\eqref{A-charact}, and for that function the above inclusion is an improvement of the inclusion $\mathcal{H}_{J,0}(\Omega)\subset L^2(\Omega)$. In fact, as we will see, we have $\mathcal{L}_{\mathcal{A},p}(\Omega)\subsetneqq L^p(\Omega)$ for every $p\ge1$.

As a precedent we have Peetre's result \cite{Peetre} which asserts that $H^{\alpha/2}_0(\Omega)$ is contained in the Lorentz space $L_{\frac{2N}{N-\alpha},2}(\Omega)$, $N>\alpha$, which corresponds to taking $\mathcal{A}(s)=s^{\frac{N-\alpha}N}$.

We also recall that, since the norm in $\mathcal{L}_{\mathcal{A},p}(\Omega)$ depends only on the distribution function, it is invariant under rearrangement. We now show that the Lorentz space $\mathcal{L}_{\mathcal{A},p}(\Omega)$ is in fact an $L^p$ space with weight when restricted to radially symmetric decreasing functions, $\mathcal{L}_{\mathcal{A},p,\text{radial}}(\Omega)=L^p_{\text{radial}}(\Omega;\psi)$, where $\psi$ and $\mathcal{A}$ are related by the formula~\eqref{A-charact}.

\begin{lem}\label{lem.charact}
 If $u$ and $\psi$  are  non-negative, radially symmetric decreasing functions with compact support, then
 \begin{equation*}\label{charact}
 \int_{\mathbb{R}^N} |u(x)|^p\psi(x)\,dx=  \|u\|_{\mathcal{A},p}^p,
 \end{equation*}
 where
\begin{equation}\label{A-charact}
\mathcal{A}(s)=\int_0^s\psi\left((z/\omega_N)^{1/N}\right)\,dz.
\end{equation}
\end{lem}
\noindent{\it Proof.}
Using  the layer cake representation \eqref{layercake} of $u^p$, we have
$$
\begin{array}{l}
\displaystyle\int_{\mathbb{R}^N} |u(x)|^p\psi(x)\,dx\displaystyle=
\int_0^\infty\int_{\mathbb{R}^N} \mathds{1}_{\{u(x)> t^{1/p}\}}\psi(x)\,dxdt \\ [4mm]
\qquad\displaystyle=\int_0^\infty N\omega_N\int_0^{(\mu(t^{1/p})/\omega_N)^{1/N}} r^{N-1}\psi(r)\,drdt=\int_0^\infty \int_0^{\mu(t^{1/p})} \psi\left((z/\omega_N)^{1/N}\right)\,dzdt\\ [4mm]
\qquad\displaystyle=p\int_0^\infty \int_0^{\mu(\tau)} \psi\left((z/\omega_N)^{1/N}\right)\,\tau^{p-1}\,dzd\tau=p\int_0^\infty \mathcal{A}(\mu(\tau))\,\tau^{p-1}\,d\tau.
\end{array}
$$
\qed

In general we have
$$
\|u\|_{\mathcal{A},p}=\|u^*\|_{\mathcal{A},p}=\|u^*\|_{L^p(\Omega;\psi)}\ge\|u\|_{L^p(\Omega;\psi)}.
$$

This characterization allows to see easily when the space $\mathcal{L}_{\mathcal{A},p}(\Omega)$ is strictly smaller than $L^p(\Omega)$.
\begin{prop}\label{lorentz-sharp}
If $\mathcal{A}'(0^+)=\infty$ and $\lim_{s\to0^+}\frac{\mathcal{A}'\mathcal{A}'''}{(\mathcal{A}'')^2}\ge\nu>1$, then $\mathcal{L}_{\mathcal{A},p}(\Omega)\varsubsetneq L^p(\Omega)$ for every $p\ge1$.
\end{prop}
\noindent{\it Proof.} The inclusion is immediate since $\mathcal{A}(s)\ge cs$, and thus
$$
\int_0^\infty \mathcal{A}(\mu(t))\,t^{p-1}\,dt\ge c\int_0^\infty \mu(t)\,t^{p-1}\,dt=c\|u\|_p^p.
$$
To get that the inclusion is proper assume for simplicity that there exists $\rho>0$ small such that $\Omega$ contains the ball $B_\rho$, and consider the function
$$
u(x)=v^{1/p}(|x|)\mathds{1}_{B_\rho}, \qquad v(s)=\frac{-\psi'(s)}{s^{N-1}\psi^\nu(s)},\quad 1<\nu<2,
$$
where $\psi(s)=\mathcal{A}'(\omega_Ns^N)$. The condition on $\mathcal{A}$ implies that if $\rho$ is small then $v$ is decreasing in $(0,\rho)$.
We first have $u\in L^p(\Omega)$,
$$
\int_{\Omega} |u(x)|^p\,dx=N\omega_N\int_0^{\rho}\frac{-\psi'(s)}{\psi^\nu(s)}\,ds=N\omega_N\int_{\psi(\rho)}^\infty\frac{dr}{r^\nu}\,dr<\infty.
$$
On the other hand, $u\notin L^p(\Omega;\psi)$, since
$$
\int_{\Omega} |u(x)|^p\psi(x)\,dx=N\omega_N\int_0^{\rho}\frac{-\psi'(s)}{\psi^{\nu-1}(s)}\,ds=
N\omega_N\int_{\psi(\rho)}^\infty\frac{dr}{r^{\nu-1}}\,dr=\infty.
$$
\qed

If $J$ satisfies hypotheses  \eqref{kernel2}--\eqref{slowly}, then it is an exercise to check that in fact $\mathcal{A}$ satisfies the hypotheses of Proposition~\ref{lorentz-sharp}. First observe that property~\eqref{slowly2} implies

$$
\lim_{s\to0^+}\frac{M(s)}{\ell(s)}=\infty.
$$
(This could be deduced directly without~\eqref{slowly2} by Karamata's Theorem, see for instance~\cite{BGT}, since $\ell$ is a slowly varying function,).
Then using again~\eqref{slowly2}
$$
\lim_{s\to0^+}\frac{\mathcal{A}'(s)\mathcal{A}'''(s)}{(\mathcal{A}''(s))^2}=
\lim_{s\to0^+}\frac{M(s)(sM''(s)-(N-1)M'(s)}{(M'(s))^2}=
\lim_{s\to0^+}\frac{M(s)}{\ell(s)}\left(N-\frac{s\ell'(s)}{\ell(s)}\right)=\infty.
$$

\noindent{\it Proof of Theorem \ref{th-Lorentz}.}

The inclusion in the Lorentz space $\mathcal{L}_{\mathcal{A},2}(\Omega)$ is obtained by Theorems~\ref{th-hardy1} and~\ref{th-symmet}, together with Lemma~\ref{lem.charact}. In fact, if $u\in\mathcal{H}_{J,0}(\Omega)$ let $u^*$ be its decreasing rearrangement, defined in $B_R$. Then
$$
\mathcal{E}(u,u)\ge \mathcal{E}_{\mathcal{K}}(u,u)\ge\mathcal{E}_{\mathcal{K}}(u^*,u^*)\ge c\int_{B_R} |u^*(x)|^2M(\rho|x|/R)\,dx=  c\|u^*\|_{\mathcal{A},2}^2
=  c\|u\|_{\mathcal{A},2}^2
$$
where $\mathcal{A}$ is defined in \eqref{A-charact} with $\psi(x)=M(\rho|x|/R)$.
\qed

\section{The nonlocal operator $\mathfrak{L}$}\label{sect-eigenvalues}
\setcounter{equation}{0}

We present in this section some properties of the linear nonlocal operator $\mathfrak{L}$ defined from the kernel $J$ by \eqref{operatorL}. Observe that, as it occurs with the bilinear form, this operator is different from the censored  nonlocal operator defined by
$$
L_cu(x)=\textrm{P.V.}
\int_{\Omega}(u(x)-u(y))J(x,y)\,dy,
$$
even if $u$ vanishes outside $\Omega$.
In fact $\mathfrak{L}u=L_cu+\Lambda u,$ whe $\Lambda$ is the function defined in~\eqref{Lambda}.

\subsection{Regularity properties of $\mathfrak{L}$}

Here we are in particular interested in operators with a kernel of the form  $J(x,y)=\mathcal{K}(x-y)$, where $\mathcal{K}$ satisfies hypotheses  \eqref{kernel2}--\eqref{slowly}, but in some cases only the lower estimate $J(x,y)\ge\mathcal{K}(x-y)$ is needed.

We first study the effect of applying the operator $\mathfrak{L}$ to a H\"older continuous function.

\begin{theorem}\label{C_norm}
	Assume $J(x,y)=\mathcal{K}(x-y)$, where $\mathcal{K}$ satisfies  \eqref{kernel2}--\eqref{slowly}. If $u\in C^{\nu}(\R^N)$ for some $\nu\,\in\, (0,1)$, then $\mathfrak{L}u\in C(\R^N)$, with some modulus of continuity that depends on $\mathcal{K}$ and $\nu$.
\end{theorem}
\noindent{\it Proof.}
Let us estimate the difference $|\mathfrak{L}u(x)-\mathfrak{L}u(y)|$ for $x,y \in \R^N$. Let $R<1$ to be fixed.

$$
\left|\mathfrak{L}u(x)-\mathfrak{L}u(y)\right|=\left|\int_{\R^N}\left(u(x)-u(x+z)-u(y)+u(y+z)\right)\mathcal{K}(z)\,dz\right|\leq I_1 +I_2,
$$

where
$$
I_1=\int_{B_R}\left(\left|u(x)-u(x+z)\right|+\left|u(y)-u(y+z)\right|\right)\mathcal{K}(z)\,dz,
$$

$$
I_2=\int_{B_R^c}\left(\left|u(x)-u(y)\right|+\left|u(x+z)-u(y+z)\right|\right)\mathcal{K}(z)\,dz.
$$

For $I_1$, using that $\left|u(x)-u(x+z)\right|\leq \left[u\right]_{C^\nu}|z|^\nu$ we get

$$
I_1\leq 2\left[u\right]_{C^\nu}\left|\int_{B_R}|z|^\nu\mathcal{K}(z)\,dz \right|=c\int_0^R\frac{\ell(s)}{s^{1-\nu}}\,ds=A(R).
$$

For $I_2$, using that $\left|u(x+z)-u(y+z)\right|\leq \left[u\right]_{C^\nu}|x-y|^\nu$ we get

$$
I_2\leq 2\left[u\right]_{C^\nu}|x-y|^\nu\int_{B_R^c}\mathcal{K}(z)\,dz \le c|x-y|^\nu M(R).
$$
Thus picking $R=g^{-1}(|x-y|)$, where $g(R)=\left(A(R)/M(R)\right)^{1/\nu}$, we obtain
$$
\left|\mathfrak{L}u(x)-\mathfrak{L}u(y)\right|\leq 2c\varpi(|x-y|)
$$
where $\varpi=M\circ g^{-1}$.
\qed

When  $\ell(0)<\infty$  we have  $A(R)\sim R^\nu$ and $M(R)\sim \log(1/R)$ for $R\to0$.  In that case $g(R)\sim R\log^{-1/\nu}(1/R)$. We have then that the regularity of $\mathfrak{L}u$ is almost the same as that of $u$; in particular $\mathfrak{L}u\in C^{\nu-\varepsilon}(\mathbb{R}^N)$ for every $0<\varepsilon<\nu$.

With the same technique we can obtain the following.
\begin{cor}\label{sharp-regularity}
  If $u$ is a continuous function with a modulus of continuity $\varpi_0$, then $\mathfrak{L}u$ is continuous provided
  $$
  \lim_{R\to0}\int_0^R\frac{\varpi_0(s)\ell(s)}s\,ds=0.
  $$
\end{cor}

Also of interest is to obtain integrability properties of $\mathfrak{L}u$ when $u$ is the characteristic function of some set, depending on the regularity of the boundary. Observe that  for $E\subset\mathbb{R}^N$
$$
\int_{E}%
\mathfrak{L}\mathds{1}_E(x)\,dx=\int_{E}%
\int_{E^c}\mathcal{K}(x-y)\,dydx =\mathcal{E}_1(\mathds{1}_E,\mathds{1}_E).
$$
If $E\subset\Omega$ this quantity  coincides with $\mathcal{E}(\mathds{1}_E,\mathds{1}_E)$.
This is called the $J$--perimeter of the set $E$, $P_J(E)$. See \cite{CaffarelliRoquejoffreSavin} for the fractional perimeter and for more general definitions when $E\not\subset\Omega$.

We say that $\partial E$ has a modulus of continuity $\varpi_0$ if it is locally  the graph of a  function  $\eta$ defined on a small ball $B\subset\mathbb{R}^{N-1}$,
such that
$$
|\eta(z_1)-\eta(z_2)|\le\varpi_0(|z_1-z_2|),\qquad\text{for all } z_1,\,z_2\in B.
$$

\begin{theorem}\label{perimeter}
  Assume $J(x,y)=\mathcal{K}(x-y)$, where $\mathcal{K}$ satisfies  \eqref{kernel2}--\eqref{slowly}. If $E$ has modulus of continuity $\varpi_0$, then $P_J(E)<\infty$  provided
 $$
  \int_{0}^\rho\frac{\varpi_0(s)\ell(s)}s\,ds<\infty.
  $$
\end{theorem}
\noindent{\it Proof.} We estimate $\int_{E^c}\mathcal{K}(x-y)\,dy$ for $x\in E$.
To that purpose let $x\in E$, with $\delta(x)=\text{dist}(x,\partial E)=r>0$, and let $B_r=\{|y-x|<r\}$. We have
$$
\int_{E^c}\mathcal{K}(x-y)\,dy\le\int_{B_r^c}\mathcal{K}(x-y)\,dy= \int_{r<|z|<\rho}\mathcal{K}(z)\,dz+\int_{|z|>\rho}\mathcal{K}(z)\,dz=c_1M(r)+c_2,
$$
if $r<\rho$; if $r\ge\rho$ we directly deduce $
\int_{E^c}\mathcal{K}(x-y)\,dy<c$.

As an immediate consequence we have that if $D_0=\{x\in E,\, \delta(x)\ge\rho\}$, then $\int_{D_0}\mathfrak{L}\mathds{1}_E\le c$.
Then it suffices to show that $\int_DM(\delta(x))\,dx<\infty$ for each set of the form
$$
D=\left\{(x',x_N)\in\mathbb{R}^{N-1}\times\mathbb{R},\,|x'|<\varepsilon,\,0< x_N<\eta(x')\right\},
$$
where $\eta<\rho$. The regularity of the function $\eta$ gives that
$$
\delta(x)\ge \varpi_0^{-1}(|\eta(x')-x_N|).
$$
Therefore, since $M$ is  decreasing,
$$
\begin{array}{rl}
\displaystyle \int_DM(\delta(x))\,dx&\displaystyle=\int_{|x'|<\varepsilon}\int_0^{\eta(x')}M(\delta(x',x_N))\,dx_Ndx' \\ [4mm]
&\displaystyle\le\int_{|x'|<\varepsilon}\int_0^{\eta(x')}M(\varpi_0^{-1}(|\eta(x')-x_N|))\,dx_Ndx'\\ [4mm]
&\displaystyle\le c\int_0^{\rho}M(\varpi_0^{-1}(z))\,dz=c\int_0^\rho\frac{\varpi_0(s)\ell(s)}s\,ds.
\end{array}
$$
\qed

Observe that if $\partial E$ is a domain with very weak continuity, say $\varpi_0(s)=(\log1/s)^{-\sigma}$, $\sigma>0$, there always exists an operator $\mathfrak{L}$ with a singularity so weak that makes $\mathfrak{L}\mathds{1}_E$ still integrable, just choose $\ell(s)=(\log1/s)^{-1}$ for that set. Therefore  the characteristic function of a bounded domain $E\subset\Omega$ with a Lebesgue spine belongs to some Sobolev type space $\mathcal{H}_{J,0}(\Omega)$, and the set has finite $J$-perimeter, where $J(x,y)=|x-y|^{-N}(\log1/|x-y|)^{-1}$.

We end this subsection estimating the action of $\mathfrak{L}$ to a specific power $|x|^{-N/2}$, precisely the one needed in the proof of Hardy inequality. 

\begin{lem}\label{lema.x-N}
  Assuming hypotheses  \eqref{kernel2}--\eqref{slowly}, there exists $\varepsilon>0$ small such that
  \begin{equation}\label{Lww}
    |x|^{N/2}\mathfrak{L}|x|^{-N/2}\ge c M(|x|)\qquad\text{for } 0<|x|<\varepsilon.
  \end{equation}
\end{lem}
\noindent{\it Proof.} By the hypotheses we have $J(x,y)\ge c\ell(|x-y|)|x-y|^{-N}$ for $|x-y|<\rho$ for some $\rho>0$. We thus have
$$
\begin{array}{rl}
|x|^{N/2}\mathfrak{L}|x|^{-N/2}&\displaystyle\ge
c\int_{|x-y|<\rho,\,|y|>|x|}\left(1-\frac{|x|^{N/2}}{|y|^{N/2}}\right)\frac{\ell(|x-y|)}{|x-y|^{N}}\,dy \\ [4mm] &\displaystyle- c\int_{|y|<|x|}\left(\frac{|x|^{N/2}}{|y|^{N/2}}-1\right)\frac{\ell(|x-y|)}{|x-y|^{N}}\,dy
\\ [4mm] &\displaystyle+
\int_{|x-y|>\rho}\left(1-\frac{|x|^{N/2}}{|y|^{N/2}}\right)J(x,y)\,dy\\ [4mm] &\displaystyle=I_1-I_2+I_3.
\end{array}
$$
We see that $I_2$ is convergent and $I_3$ is positive if for instance $|x|<\rho/2$. As to $I_1$  we get
$$
\begin{array}{rl}
I_1&\displaystyle=c\int_{|z|<\rho,\,|z-x|>|x|}
\frac{|z-x|^{N/2}-|x|^{N/2}}{|z-x|^{N/2}}\,\frac{\ell(|z|)}{|z|^{N}}\,dz\\ [4mm]
&\displaystyle\ge
c\int_{3|x|<|z|<\rho}\frac{\ell(|z|)}{|z|^{N}}\,dz\ge cM(3|x|)\ge cM(|x|),
\end{array}$$
if $|x|<\rho/3$, where we have used in the last inequality that $\ell$ is slowly varying. We conclude since $M(0^+)=\infty$.
\qed

In the fractional Laplacian case it is easy to obtain, by means of the Fourier transform,
$$
(-\Delta)^{\alpha/2}|x|^\gamma=c_{N,\alpha,\gamma}|x|^{\gamma-\alpha},
$$
whenever $\gamma>-N$. In particular
$$
|x|^{\frac{N-\alpha}2}(-\Delta)^{\alpha/2}|x|^{-\frac{N-\alpha}2}=c_{N,\alpha}|x|^{-\alpha},
$$
which gives the weight for the fractional Hardy inequality~\eqref{hardy-lapfrac}. Also the sharp constant can be obtained from that identity. Compared with estimate~\eqref{Lww} we formally try to put $\alpha=0$ on the left, taking care of the constant, obtaining a logarithmic type function on the right.

\subsection{The eigenvalue problem}

We now  consider the problem of finding the eigenvalues and eigenfunctions of $\mathfrak{L}$ in $\mathcal{H}_{J,0}(\Omega)$, that is
\begin{equation}\label{eigen}
\begin{cases}
\mathfrak{L}\varphi=\lambda\varphi,& \mbox{in }\Omega,\\ \varphi=0,& \mbox{in } \Omega^c.
\end{cases}
\end{equation}
We only need to consider here the lower bound $J(x,y)\ge \mathcal{K}(x-y)$.

\begin{theorem}\label{eigen2}
Assume $J$ satisfies hypotheses \eqref{kernel2}--\eqref{slowly}. The first eigenvalue
$$
  \lambda_1=\min_{\substack{\varphi\in \mathcal{H}_{J,0}(\Omega)\\ \|\varphi\|_2=1} }\mathcal{E}(\varphi,\varphi),
  $$
  is positive and isolated, and the first eigenfunction does not vanish in $\Omega$.
 \end{theorem}

The proof is rather standard so we only sketch the main steps.

\noindent{\it Proof.}
Consider the functional $\Psi\,:\,\mathcal{H}_{J,0}(\Omega)\to\mathbb{R}^+$, defined by $\Psi(u)=\mathcal{E}(u,u)$,
and let $\mathcal{M}=\{u\in\mathcal{H}_{J,0}(\Omega)\,:\,\|u\|_2=1\}$. Let $\{u_k\}$ be a minimizing sequence for $\Psi$ in $\mathcal{M}$, that is
$$
\lim_{k\to\infty}\Psi(u_k)=c=\inf_{u\in\mathcal{M}}\Psi(u)\ge0.
$$
Then $\{u_k\}$ is bounded in $\mathcal{H}_{J,0}(\Omega)$, so there exists a subsequence (still denoted by $\{u_k\}$), such that $u_k\rightharpoonup u^*$ in $\mathcal{H}_{J,0}(\Omega)$, and also
$$
\lim_{k\to\infty}\mathcal{E}(u_k,\eta) =\mathcal{E}(u^*,\eta),\qquad\text{for every } \eta\in\mathcal{H}_{J,0}(\Omega).
$$
By Theorem \ref{th-compact}, there exists a new subsequence converging to $u^*$ in $L^2(\Omega)$, so $\|u^*\|_2=1$ and $u^*\in\mathcal{M}$. This gives
$$
c=\lim_{j\to\infty}\Psi(u_j)\ge\Psi(u^*)\ge c,
$$
and $\Psi(u^*)= c $. The first eigenvalue is then $\lambda_1=\Psi(u^*)>0$, with corresponding eigenfunction $\varphi_1=u^*$. The fact that $u^*\ge0$ or $u^*\le0$ follows from~\eqref{absolute}. Regularity of the first eigenfucntion, which is obtained in the next section, would in fact imply that $u^*$ does not vanish in $\Omega$. Observe that if $u^*(x)\ge0$ and $u^*(x_0)=0$, then
$$
0=\lambda_1u^*(x_0)=-\int_{\mathbb{R}^N}u^*(y)J(x,y)\,dy<0.
$$
Finally suppose that there exists $v\in \mathcal{H}_{J,0}(\Omega)$ with $\|v\|_2=1$ such that $\mathfrak{L}v=\lambda_1v$. Then $w=v-u^*$ also satisfies $\mathfrak{L}w=\lambda_1w$, and thus it has a definite sign. This gives $|v|\ge |u^*|$ or the opposite $|v|\le |u^*|$. But they have equal $L^2$ norm, so $|v|=|u^*|$, and thus $v=\pm u^*$, that is, $\lambda_1$ is isolated.
\qed

We also have
\begin{theorem}\label{prop-eigen}
In the above hypotheses there exists a sequence of eigenvalues $\{\lambda_j\}$ and eigenfunctions $\{\varphi_j\}$ to problem~\eqref{eigen} with the following properties:
\begin{enumerate}
  \item $\{\lambda_j\}$ is nondecreasing with limit $\infty$.
  \item If $P_j=\{\varphi\in \mathcal{H}_{J,0}(\Omega),\,\varphi\neq0,\;\mathcal{E}(\varphi,\varphi_k)=0\;\forall\;k=1,\cdots,j-1\}$, then
  $$
  \lambda_j=\min_{ \stackrel{\varphi\in P_j}{\|\varphi\|_2=1}}\mathcal{E}(\varphi,\varphi)=\mathcal{E}(\varphi_j,\varphi_j).
  $$
  \item The eigenfunctions form an orthonormal basis of $\mathcal{H}_{J,0}(\Omega)$, that is, for every $u\in\mathcal{H}_{J,0}(\Omega)$, it holds
  $$
  \lim_{n\to\infty}\mathcal{E}(u-\sum_{j=1}^n\verb"u"_j\varphi_j,\eta)=0,\qquad\text{for every } \eta\in\mathcal{H}_{J,0}(\Omega),
  $$
  where
  $$
  \verb"u"_j=\frac{1}{\lambda_j}\mathcal{E}(u,\varphi_j).
  $$
\end{enumerate}
\end{theorem}
\noindent{\it Proof.}
The same construction as before gives the existence of the sequence of eigenvalues and eigenfunctions, with $\lambda_j=\Psi(\varphi_j)$. Two eigenfunctions $\varphi,\,\psi$, corresponding to two different eigenvalues $\lambda,\,\mu$ are orthogonal, in  $L^2(\Omega)$ and $\mathcal{H}_{J,0}(\Omega)$ since
$$
\lambda\int_{\Omega}\varphi\psi=\mathcal{E}(\varphi,\psi)=\mu\int_{\Omega}\varphi\psi,
$$
so that $\int_{\Omega}\varphi\psi=0$ and as a consequence $\mathcal{E}(\varphi,\psi)=0$.

If the sequence $\{\lambda_j\}$ were bounded there would exist a subsequence of $\{\varphi_j\}$ converging in $L^2(\Omega)$, but orthogonality in $L^2(\Omega)$ implies
$\|\varphi_j-\varphi_k\|_2=2$, for every $j,\,k$, which is a contradiction.
\qed

  On the other hand, we can describe the space $\mathcal{H}_{J,0}(\Omega)$,  the operator $\mathfrak{L}$ and the bilinear form $\mathcal{E}$ in terms of the eigenvalues.
\begin{prop}\label{eigen-op} In the above hypotheses,
  $$
  \mathcal{H}_{J,0}(\Omega)=\left\{u\in L^2(\Omega),\, u\equiv0 \text{ in } \Omega^c,\; \|u\|_{\mathcal{H}_{J,0}}\equiv\left(\sum_{j=1}^\infty\lambda_j\verb"u"_j^2\right)^{1/2}<\infty\right\}.
  $$
  and
  $$
  \begin{array}{l}
\displaystyle  \mathfrak{L}u=\sum_{j=1}^\infty\lambda_j\verb"u"_j\varphi_j,\qquad \mathfrak{L}^{1/2}u=\sum_{j=1}^\infty\lambda_j^{1/2}\verb"u"_j\varphi_j,\\ [4mm] \displaystyle\mathcal{E}(u,v)=\sum_{j=1}^\infty\lambda_j\verb"u"_j\verb"v"_j=\int_{\mathbb{R}^N}\mathfrak{L}^{1/2} u\,\mathfrak{L}^{1/2}v\,,
\end{array}
  $$
  where
  $$
  \verb"u"_j=\int_\Omega u\varphi_j=\frac{1}{\lambda_j}\mathcal{E}(u,\varphi_j),\quad
  \verb"v"_j=\int_\Omega v\varphi_j=\frac{1}{\lambda_j}\mathcal{E}(v,\varphi_j).
  $$
\end{prop}

Observe that $\|u\|_{\mathcal{H}_{J,0}}=\|\mathfrak{L}^{1/2}u\|_{L^2}$. With this construction we have that $\mathfrak{L}\,:\,\mathcal{H}_{J,0}(\Omega)\to H^*(\Omega)$ is an isomorfism, where $H^*(\Omega)$ is the closure of the set of functions $v=\sum_{j=1}^\infty \verb"u"_j\varphi_j$ with the norm $\|v\|_{H^*}=\left(\sum_{j=1}^\infty\lambda_j^{-1}\verb"v"_j^2\right)^{1/2}$. The duality product is
$$\langle u,v\rangle_{\mathcal{H}_{J,0}\times H^*}=\sum_{j=1}^\infty\lambda_j\verb"u"_j\lambda_j^{-1}\verb"v"_j=\int_\Omega uv.$$

We finally estimate the first eigenvalue in terms of the size of the domain. We obtain a lower bound of the type of the one obtained in \cite{Li-Yau} in the case of the Laplacian and in \cite{Yolcu-Yolcu} for the fractional Laplacian case. We  consider the multiplier $m(\xi)$ associated to the kernel $\mathcal{K}$, see Subsection~\ref{sect-inclusion}. Observe that $m$ is radial, $m(\xi)=\mathfrak{m}(|\xi|)$, and if $\ell$ is nonincreasing then  $\mathfrak{m}$ is nondecreasing. To see that, for each $\xi\in\mathbb{R}^N$ given, let $R=R_\xi$ be  the rotation that carries $\xi$ to the first  axis, that is, $R(\xi)=|\xi|e_1$, where $e_1=(1,0,\dots,0)$, and let $S=R^{-1}$. If we put $z=S^Ty|\xi|$, then
$$
y\cdot \xi=y\cdot S(|\xi|e_1)=(S^Ty)\cdot( |\xi|e_1)=z_1.
$$
Thus
$$
m(\xi)=\int_{\mathbb{R}^N}\left(1-\cos (y\cdot \xi)\right)\mathcal{K}(y)\,dy=\int_{\mathbb{R}^N}\frac{1-\cos z_1}{|z|^{N}}\ell(|\xi|^{-1}|z|)\,dz.
$$
Thus $\mathfrak{m}$ increases when $\ell$ decreases. Put now
$$
g(t)=\int_{|\xi|\le t}m(\xi)\,d\xi.
$$
We need to suppose that $g$ satisfies
\begin{equation}\label{g-decreasing}
Ng(t)\le tg'(t).
\end{equation}
In particular it implies
$$
\begin{cases}
  \mathcal{K}(z)\ge c|z|^{-N}&\text{ if } |z|\le1,\\
  \mathcal{K}(z)\le c|z|^{-N}&\text{ if } |z|\ge1,
\end{cases}
$$
which is not too restrictive in our situation of nonintegrable L\'evy kernels.

\begin{theorem}\label{berezin}
 Assume $J$ satisfies hypotheses \eqref{kernel2}--\eqref{slowly}, with $\ell$ nonincreasing and $g$ satisfies~\eqref{g-decreasing}.  Then
 \begin{equation}\label{dec-eigen}
  \lambda_1\ge\frac{|\Omega|}{(2\pi)^N}g\left(\frac{2\pi}{\left(\omega_N|\Omega|\right)^{1/N}}\right).
 \end{equation}
\end{theorem}
\noindent{\it Proof.}
The first eigenvalue satisfies
$$
\lambda_1=\mathcal{E}(\varphi_1,\varphi_1)\ge\mathcal{E}_{\mathcal{K}}(\varphi_1,\varphi_1)=
\int_{\mathbb{R}^N}m(\xi)|\widehat\varphi_1(\xi)|^2\,d\xi,
$$
where $\|\varphi_1\|_2=\|\widehat\varphi_1\|_2=1$.
Put $h(\xi)=\|\widehat\varphi_1\|_\infty^2\mathds{1}_{\{|\xi|<K\}}$, for some $K$ to be determined. Since $\mathfrak{m}$ is increasing, we have
$$
\left(\mathfrak{m}(|\xi|)-\mathfrak{m}(K)\right)\left(|\widehat\varphi_1(\xi)|^2-h(\xi)\right)\ge0.
$$
Therefore
$$
\mathfrak{m}(K)\left(|\widehat\varphi_1(\xi)|^2-h(\xi)\right)
\le\mathfrak{m}(|\xi|)\left(|\widehat\varphi_1(\xi)|^2-h(\xi)\right).
$$
Integrating this inequality we get
$$
\int_{\mathbb{K}^N}\mathfrak{m}(R)\left(|\widehat\varphi_1(\xi)|^2-h(\xi)\right)\,d\xi
\le\int_{\mathbb{R}^N}\mathfrak{m}(|\xi|)\left(|\widehat\varphi_1(\xi)|^2-h(\xi)\right)\,d\xi\le0,
$$
provided $K$ is chosen such  that $g(K)= \frac{\lambda_1}{\|\widehat\varphi_1\|_\infty^2}$. We get
$$
1=\int_{\mathbb{R}^N}|\widehat\varphi_1(\xi)|^2\,d\xi\le\int_{\mathbb{R}^N}h(\xi)\,d\xi=
\|\widehat\varphi_1\|_\infty^2K^N\omega_N,
$$
and therefore
$$
\lambda_1\ge\|\widehat\varphi_1\|_\infty^2g\left(\left(\omega_N\|\widehat\varphi_1\|_\infty^2\right)^{-1/N}\right).
$$
The function on the right is decreasing by~\eqref{g-decreasing},  and since $\varphi_1$ hast compact support contained in $\Omega$, by Cauchy-Schwartz inequality
$$
\|\widehat\varphi_1\|_\infty^2\le \frac{|\Omega|}{(2\pi)^N},
$$
so that we conclude \eqref{dec-eigen}.
\qed

{\bf Remarks.} $i)$ Since always $m(\xi)\ge c|\xi|^2$ near the origin, then \eqref{dec-eigen} means, for large domains,
$$
\lambda_1\ge c|\Omega|^{-\frac{2}{N}}.
$$
If moreover $\mathcal{K}(z)\ge c|z|^{-\alpha}$ for some $\alpha>0$, then
$$
\lambda_1\ge c|\Omega|^{-\frac{\min\{\alpha,2\}}{N}}.
$$
$ii)$ With the same technique it can also be obtained the estimate for the sum of the eigenvalues
$$
 \sum_{j=1}^k\lambda_j\ge\frac{|\Omega|}{(2\pi)^N}g\left(\frac{2\pi k^{1/n}}{\left(\omega_N|\Omega|\right)^{1/N}}\right).
$$

\section{Elliptic problems}\label{sect-problems}
\setcounter{equation}{0}

In this section we explain some results on integral regularity of solutions to elliptic problems of the form
\begin{equation}\label{problem}
\begin{cases}
\mathfrak{L}u=f,& \mbox{in }\Omega,\\ u=0,& \mbox{in }\Omega^c,
\end{cases}
\end{equation}
Here $\mathfrak{L}$ is the operator \eqref{operatorL}. Given a function $f\in H^*(\Omega)$, we say that $u:\R^N\to \R$ is a weak solution of \eqref{problem}, if $u\in \mathcal{H}_{J,0}(\Omega)$ is a function such that
\begin{equation}\label{weaksol-D0}
\mathcal{E}(u,\phi)=\int_{\Omega}f\phi,\quad\mbox{for all } \phi \in \mathcal{H}_{J,0}(\Omega).
\end{equation}
Existence and uniqueness of solution is proved in \cite{Felsinger-Kassmann-Voigt} in a more general framework. In fact we only need a Poincar\'e inequality, and then the proof is standard. We include it for completeness.
\begin{prop}\label{prop-lineal}
	Assume Poincar\'e inequality \eqref{poincare} holds. Then problem \eqref{problem} has a unique weak solution $u\,\in\,\mathcal{H}_{J,0}(\Omega)$.
\end{prop}
\noindent{\it Proof.}
Consider the energy functional $\mathcal{G}\,:\,\mathcal{H}_{J,0}(\Omega)\to\mathbb{R}$ associated to the problem \eqref{problem}, defined by
\begin{equation*}\label{energy}
\mathcal{G}(u)=\frac{1}{2}\mathcal{E}(u,u)-\int_{\Omega}fu.
\end{equation*}
This functional  is well defined thanks to Poincar\'e inequality, it is Fr\'echet differentiable in $u\,\in\, \mathcal{H}_{J,0}(\Omega)$ and for any $\phi \in\mathcal{H}_{J,0}(\Omega)$
$$
\langle \mathcal{G}'(u),\phi\rangle=\mathcal{E}(u,\phi)-\int_{\Omega} f\phi,
$$
that is, critical points  of $\mathcal{G}$ are weak solutions to  \eqref{problem}. The result is obtained by minimizing the functional $\mathcal{G}$. Observe also that $\mathcal{H}_{J,0}(\Omega)$ is a Hilbert space so we could have used Riesz representation theorem.
\qed

Maximum principle and comparison principle for weak solutions (or more generally for supersolutions) to \eqref{problem} are also easy to obtain.

 \begin{prop}\label{maximum}
 	If $u\in\mathcal{H}_J(\mathbb{R}^N)$ then $\mathfrak{L}u\geq 0$ in $\Omega$ and $u\ge0$ in $\Omega^c$ imply $u\geq0$ in $\Omega$.
 \end{prop}
\noindent{\it Proof.} Property $\mathfrak{L}u\geq 0$ in $\Omega$ actually means $\mathcal{E}(u, \phi)\ge0$ for every $\phi \in \mathcal{H}_J(\Omega)$, $\phi\ge0$. Since $u^-\ge0$ and $u^-\in\mathcal{H}_J(\Omega)$, by the Stroock-Varopoulos inequality we have
$$
0\ge- \mathcal{E}(u^-, u^-)\ge \mathcal{E}(u, u^-)\ge 0.
$$
Hence $u^-\equiv0$.
\qed

The comparison principle follows immediately as a consequence.

The following result, due to Kassmann and Mimica  \cite{Kassmann-Mimica}, explain the weak character of the smoothing effect in problem~\eqref{problem}. This result is sharp by Corollary~\ref{sharp-regularity}
\begin{theorem}
	Assume hypotheses \eqref{kernel2}--\eqref{slowly} and let $u$ be a bounded weak solution to \eqref{problem} with $f\in L^{\infty}(\Omega)$. Then there exist constants $c>0$ and $\beta \in (0,1)$ such that
	\begin{equation*}
	|u(x)-u(y)|\leq c\left(\|u\|_{\infty} + \|f\|_{\infty}\right)\varpi(|x-y|),\quad\text{for every } x,\,y\in\Omega,
	\end{equation*}
where $\varpi(s)=\frac{1}{M^\beta(s)}$.
\end{theorem}

We now study the smoothing effect in terms of integrability. Before that we show first that the solution is not worse that the datum.
In the local case $-\Delta u=f$ there is a strong smoothing effect: $u$ is bounded provided $f\in L^p(\Omega)$, $p>N/2$, see for instance \cite[Theorem 8.15]{Gilbarg-Trudinger}, from where some ideas are borrowed below. In fact, the same calculation allows to get easily the  conclusion in the fractional Laplacian framework, $(-\Delta)^{\alpha/2} u=f$, when $p>N/\alpha$. Recall that here we are in the borderline $\alpha\sim0$. It would be interesting to obtain $u\in L^\infty(\Omega)$ for $f\in L^p(\Omega)$ for every $p<\infty$, but $f\notin L^\infty(\Omega)$.

\begin{theorem}\label{Lp_bound}
	Assume hypotheses \eqref{kernel2}--\eqref{slowly} and let $u$ be a weak solution to \eqref{problem} with $f\in L^{p}(\Omega)$,  $2\leq p \leq \infty$. Then,
\begin{equation*}\label{eq:worse}
	\| u\|_{p}\leq C\| f\|_{p},
\end{equation*}
where the constant $C$ depends only on the kernel and $\Omega$.
\end{theorem}
\noindent{\it Proof.}
Consider first the case $p=\infty$.  Let $B$ be any large  ball such that $\overline\Omega  \subset B$, and let $\eta\in C^{\infty}_c(B)$ be such that, $0 \leq \eta (x) \leq 1$, $x\in\R^N$ and $\eta \equiv 1$, in $\Omega$. Then, for each $x\in \Omega$, we have

$$
\mathfrak{L}\eta(x)=\int_{\R^N} \left(\eta(x)-\eta(y)\right)J(x-y)\,dy\geq\displaystyle\int_{B^c} J(x-y)\,dy=c>0.
$$
Taking $\omega(x)=\frac{\|f\|_{\infty}}{c}\eta(x)$ we  have  $\mathfrak{L}u\leq \mathfrak{L}\omega$ in $\Omega$, and $\omega\ge0$ in $\Omega^c$. Thus by the comparison principle  we get $u\leq\omega$ in $\Omega$, and hence $u\leq C\|f\|_\infty$. Similarly we have that $-u\leq C\|f\|_\infty$.

For the case $2\leq p <\infty$ let us see first the formal calculus. Choosing as test function $\phi=|u|^{p-2}u$, and using Poincar\'e, Stroock-Varopoulos and H\"older inequalities, we get, modulo multiplicative constants,
$$
\|u\|_p^p=\||u|^{\frac p2}\|_2\le\mathcal{E}(|u|^{\frac p2},|u|^{\frac p2})\le\mathcal{E}(u,|u|^{p-2}u)=\int f |u|^{p-2}u\le\|f\|_p\|u\|_p^{p-1}.
$$
We would get the result if $\|u\|_p$ is finite. Also $\phi=|u|^{p-2}u$ is not an admissible test function. The justification works as usual through truncation, see for instance \cite{Gilbarg-Trudinger}.
Let us consider for any $T>0$ the  function
\begin{equation}
F(s)=F_T(s)=
\begin{cases}
|s|^{\frac p2}& \mbox{if }|s|\leq T ,\\ \frac p2 T^{\frac p2-1}(|s|-T)+T^{\frac p2}& \mbox{if }|s|>T.
\end{cases}
\end{equation}
Since $F$ is a Lipschitz convex function and $F(0)=0$, we have  $F(u)\in \mathcal{H}_J(\Omega)$. If we define  $G=(F^2)'$ then   $G'\geq 2(F')^2$, and hence Poincar\'e and Stroock-Varopoulos inequalities give

$$
\|F(u)\|_2^2\leq c\mathcal{E}(F(u),F(u))\leq c\mathcal{E}(u,G(u))=c\int_{\Omega}f(x)G(u(x))\,dx.
$$
Now observe that $|G(u)|\le  pF(u)^{\frac{2(p-1)}p}$, and $|G(u)|\le c|u|$ for $|u|>T$, so that $G(u)\in L^{\frac p{p-1}}(\Omega)$. Applying then H\"older inequality to the last integral we get
$$
\int_{\Omega}fG(u)\leq c \|f\|_{p}\|F(u)\|_2^{\frac{2(p-1)}p},
$$
and hence
$$
\|F(u)\|_2^{\frac 2p}\leq c \|f\|_{p}\,,
$$
with $c$ independent of $T$.
We conclude taking the limit as $T\to \infty$, since $\|F_T(u)\|_2^{\frac 2p}\to\|u\|_p$.
\qed

\begin{theorem}
  \label{th-Lp-smoothing}
  If $f\in L^p(\Omega)$, $p\ge2$, then $u\in \mathcal{L}_{\mathcal{A},p}(\Omega)$.
  \end{theorem}
\noindent{\it Proof.}
By the above proof we know that $|u|^{\frac p2}\in \mathcal{H}_J(\Omega)$, and also that
$$
\mathcal{E}(|u|^{\frac p2},|u|^{\frac p2})\le\|f\|_p\|u\|_p^{p-1}\le c\|f\|_p^p.
$$
Therefore, using Theorem \ref{th-Lorentz} and estimate \eqref{eq:worse}, we get
$$
\|u\|_{\mathcal{A},p}^p=\||u|^{\frac p2}\|_{\mathcal{A},2}^2\le c\mathcal{E}(|u|^{\frac p2},|u|^{\frac p2})\le c\|f\|_p^p.
$$
\qed

\subsection{The problem with nonhomogeneous exterior datum}\label{sub-nonhom}
Now we want to study the problem
\begin{equation}\label{problem-D1}
\begin{cases}
\mathfrak{L}u=f,& \mbox{in }\Omega,\\ u=g,& \mbox{in }\Omega^c,
\end{cases}
\end{equation}
where $f\in H^*(\Omega)$, and $g\in\mathcal{H}_J(\mathbb{R}^N)$.

We observe that when multiplying the equation by a test function, we get
$$
\int_\Omega f\varphi=\int_\Omega\int_{\mathbb{R}^N}(u(x)-u(y))\varphi(x)J(x,y)\,dydx.
$$
Since $u$ does not necessarily vanish outside $\Omega$, the right-hand side is different from $\mathcal{E}_1(u,\varphi)$, and this is the reason of the introduction of the bilinear form $\mathcal{E}$ in~\eqref{bilinear-form}.

The solution to problem \eqref{problem-D1} is a function $u\in{\mathcal{H}}_J(\Omega)$ such that $u-g\in \mathcal{H}_{J,0}(\Omega)$ and \eqref{weaksol-D0} holds. We can solve \eqref{problem-D1} by considering the problem satisfied by $w=u-g$, and noting that $\mathfrak{L}g\in H^*(\Omega)$. We remark the recent work \cite{Dyda-Kassmann}, where conditions on the data $g$ defined only on $\Omega^c$ are imposed to guarantee that the problem is well posed, i.e., $g$ can be extended properly into $\Omega$.

\section{Nonlinear elliptic problems}\label{sect-nonlinear}
\setcounter{equation}{0}
We study in this section  the nonlinear elliptic type problem
 \begin{equation}\label{sublinear}
	\begin{cases}
		\mathfrak{L}u=f(u),& \mbox{in }\Omega,\\ u> 0,& \mbox{in }\Omega,\\ u=0,& \mbox{in }\Omega^c,
	\end{cases}
\end{equation}
We first show existence in the sublinear case, i.e.,
when $f:\left[0,\infty\right) \to \R$ is a continuous function satisfying
\begin{equation}
\label{sub1}
\frac{f(t)}t \text{ is nonincreasing on } (0,\infty)\quad  \text{ and } \lim_{t\to\infty}\frac{f(t)}t=0.
\end{equation}
See \cite{BrezisOswald} for the classical case when $\mathfrak{L}=-\Delta$.

\begin{theorem}\label{th-sublinear}
Under the assumption \eqref{sub1} problem \eqref{sublinear} admits a unique solution. 	
\end{theorem}
\noindent{\it Proof.}
Consider the energy functional $\Phi:\mathcal{H}_{J,0}(\Omega)\to \R $ defined by
$$
\Phi(u)=\frac{1}{2}\mathcal{E}(u,u)-\int_{\Omega}F(u),
$$
where $F(u)=\int_{0}^{u}f.$ From \eqref{sub1} it follows that there exist $\sigma \in (0,1)$ and $a,b>0$ such that
\begin{equation*}\label{acota}
|f(t)|\leq a+b\,t^\sigma,\quad \forall t\geq 0.
\end{equation*}
This functional is well defined since $|F(u)|\leq a_1+b_1 u^{\sigma + 1}$, and then
$$
\left|\int_{\Omega}F(u)\,dx\right|\leq  C_1+C_2\|u\|_{\sigma+1}^{\sigma +1}<\infty.
$$
On the other hand, this estimate also gives coercivity since $\sigma+1<2$,
$$
\Phi(u)\geq \frac12\|u\|_{\mathcal{H}_J}^2-C_2\|u\|_{\sigma+1}^{\sigma +1}-C_1\ge \frac12\|u\|_{\mathcal{H}_J}^2-C_3\|u\|_{\mathcal{H}_J}^{\sigma +1}-C_1.
$$

Let now $\{u_k\}\subset\mathcal{H}_{J,0}(\Omega)$ be a minimizing sequence for $\Phi$; this sequence is bounded in $\mathcal{H}_{J,0}(\Omega)$, and therefore we can assume that there is a subsequence, still denoted $\{u_k\}$, such that $u_k\rightharpoonup u$ in $\mathcal{H}_{J,0}(\Omega)$, and therefore $u_k\to u$ in $L^2(\Omega)$. We thus deduce
$$
\int_{\Omega}F(u_k)\, dx \to \int_{\Omega}F(u)\, dx,
$$
so that
$$
\Phi(u)\leq \liminf_{k\to \infty}\left(\frac{1}{2}\mathcal{E}(u_k,u_k)-\int_{\Omega}F(u_k)\right)=\liminf_{k\to \infty}\Phi(u_k)=\inf_{u\in \mathcal{H}_{J,0}(\Omega)}\Phi(u).
$$
This shows that $u$ is a global minimum for $\Phi$, and hence it is a critical point, namely a solution to \eqref{sublinear}.

Uniqueness  follows a standard argument. Suppose $u_1$ and $u_2$ are two solutions of \eqref{sublinear}, and use $\varphi_1=\dfrac{u_1^2-u_2^2}{u_1}$ and $\varphi_2=\dfrac{u_1^2-u_2^2}{u_2}$ as test functions, respectively. Then
$$
\mathcal{E}(u_1,\varphi_1)-\mathcal{E}(u_2,\varphi_2)=
\int_\Omega\left(\frac{f(u_1)}{u_1}-\frac{f(u_2)}{u_2}\right)\left(u_1^2-u_2^2\right)\leq 0,
$$
 since $f(t)/t$ is nonincreasing. On the other hand, as in the proof of Theorem~\ref{th-hardy1}, we have 
$$
\mathcal{E}(u_1,\varphi_1)-\mathcal{E}(u_2,\varphi_2)=
\mathcal{E}(u_1,u_1)-\mathcal{E}\left(u_2,\frac{u_1^2}{u_2}\right)+
\mathcal{E}(u_2,u_2)-\mathcal{E}\left(u_1,\frac{u_2^2}{u_1}\right)
\geq 0.
$$
We conclude $u_1=u_2$. \qed

We now show nonexistence for supercritical reactions $f(u)$ when $\Omega$ is star-shaped, where supercritical means above some exponent  depending on the kernel. In the fractional Laplacian case the critical exponent is $p_*=\frac{N+\alpha}{N-\alpha}$, and is proved in~\cite{RosOtonSerra} by means of a Pohozaev inequality. We follow their proof and  establish an inequality adapted to our bilinear form $\mathcal{E}$.
Let, for $\lambda>1$,
\begin{equation*}
  \label{gamma}
  \gamma(\lambda)=\lambda^{-N}\sup_{\substack{x,y\in\mathbb{R}^N\\ x\ne y}}\frac{J(x/\lambda,y/\lambda)}{J(x,y)},
\end{equation*}
and assume $\gamma(\lambda)<\infty$ for $\lambda$ close to 1.
\begin{theorem}
\label{teo-pohozaev}
If $u$ is a solution to problem \eqref{sublinear} and $\Omega$ is star-shaped, then
\begin{equation}
  \label{pohozaev}
  \int_\Omega uf(u)\le\frac{2N}{N-\sigma}\int_\Omega F(u),
\end{equation}
where $\sigma=\gamma'(1^+)$ and $F'=f$.
\end{theorem}

\begin{cor}\label{cor-p^*}
  Problem \eqref{sublinear} with $f(u)=u^p$ and $\Omega$ star-shaped has no solution for any exponent $p>p_*=\frac{N+\sigma}{N-\sigma}$.
\end{cor}

In the power case (fractional Laplacian type)
$$
J(x,y)=\begin{cases}
  |x-y|^{-N-\alpha_1}&\text{ if } |x-y|<1,\\
  |x-y|^{-N-\alpha_2}&\text{ if } |x-y|>1,
\end{cases}
$$
$\alpha_1<2$, $\alpha_2>0$, we get $\sigma=\max\{\alpha_1,\alpha_2\}$.

\noindent{\it Proof of Theorem \ref{teo-pohozaev}.} We put $\phi=u_\lambda$ as test in \eqref{weaksol-D0}, where $u_\lambda(x)=u(\lambda x)$. Since $\Omega$ is star-shaped, when $\lambda>1$ we have that $u_\lambda$ vanishes outside $\Omega$, and then $u_\lambda\in \mathcal{H}_{J,0}(\Omega)$. We have then
\begin{equation}\label{weak-poho}
\mathcal{E}(u,u_\lambda)=\int_{\Omega}f(u)u_\lambda,\quad\mbox{for all } \lambda>1.
\end{equation}
We observe that, with the above definition of $\gamma(\lambda)$, we have
$$
\begin{array}{rl}
\displaystyle\mathcal{E}(u_\lambda,u_\lambda)&\displaystyle=\frac12\iint_{\mathbb{R}^{2N}}|u(\lambda x)-u(\lambda y)|^2J(x,y)\,dxdy \\ [3mm]
&\displaystyle=\frac12\lambda^{-2N}\iint_{\mathbb{R}^{2N}}|u(x)-u(y)|^2J(x/\lambda,y/\lambda)\,dxdy \\ [3mm]
&\displaystyle\le\frac12\lambda^{-N}\gamma(\lambda)\iint_{\mathbb{R}^{2N}}|u(x)-u(y)|^2J(x,y)\,dxdy \\ [3mm]
&\displaystyle=\lambda^{-N}\gamma(\lambda)\mathcal{E}(u,u),
\end{array}
$$
so that
$$
\mathcal{E}(u,u_\lambda)\le\left(\mathcal{E}(u_\lambda,u_\lambda)\right)^{1/2}\left(\mathcal{E}(u,u)\right)^{1/2}\le
\lambda^{-N/2}\sqrt{\gamma(\lambda)}\mathcal{E}(u,u).
$$
Therefore, if $I(\lambda)=\frac{\lambda^{N/2}}{\sqrt{\gamma(\lambda)}}\mathcal{E}(u,u_\lambda)$, we deduce that $I(\lambda)\le I(1)$ for $\lambda>1$, and thus $I'(1^+)\le0$.

With this information we differentiate both sides of inequality \eqref{weak-poho} with respect to $\lambda$ at $\lambda=1$. On one hand
$$
\begin{array}{rl}
\displaystyle\left.\frac d{d\lambda}\right|_{\lambda=1^+}\mathcal{E}(u,u_\lambda)&\displaystyle=\left.\frac d{d\lambda}\right|_{\lambda=1^+}\left(\lambda^{-N/2}\sqrt{\gamma(\lambda)}I(\lambda)\right) \\ [3mm]
&\displaystyle=\left(-\frac N2+\frac{\gamma'(1^+)}2\right)I(1^+)+I'(1^+)\\ [3mm]
&\displaystyle\le-\frac12(N-\gamma'(1^+))\mathcal{E}(u,u)\\ [3mm]
&\displaystyle=-\frac12(N-\gamma'(1^+))\int_\Omega f(u)u.
\end{array}
$$
On the other hand,
$$
\left.\frac d{d\lambda}\right|_{\lambda=1^+}\int_{\Omega}f(u)u_\lambda=\int_{\Omega}x\cdot\nabla u\,f(u)=-N\int_{\Omega}F(u).
$$
Putting together this two estimates we get \eqref{pohozaev}.
\qed

\section{Neumann problems}\label{sect-neumann}
\setcounter{equation}{0}

We consider in this section Neumann type problems associated to the operator \eqref{operatorL}, following the construction made in~\cite{diPierro-RosOton-Valdinoci} for the fractional Laplacian. We therefore study the problem
\begin{equation}\label{problem-neumann}
\begin{cases}
\mathfrak{L}u=f,& \mbox{in }\Omega,\\ \mathcal{N}u=0,& \mbox{in }\Omega^c,
\end{cases}
\end{equation}
where
\begin{equation*}
\label{operatorL2}\mathfrak{L}u(x)=
\int_{\mathbb{R}^{N}}(u(x)-u(y))J(x,y)\,dy,\qquad x\in\Omega,
\end{equation*}
and
\begin{equation*}
\label{operatorN}\mathcal{N}u(x)=
\int_{\Omega}(u(x)-u(y))J(x,y)\,dy,\qquad x\in\Omega^c.
\end{equation*}
The introduction of the exterior operator $\mathcal{N}$, which plays the role of a Neumann operator, is motivated by the following property, which can be interpreted as an integration by parts formula: For every $u,\,v\in{\mathcal{H}}_J(\Omega)$ it holds
\begin{equation*}
  \label{porpartes}
  \int_\Omega v\mathfrak{L}u+\int_{\Omega^c}v\mathcal{N}u={\mathcal{E}}(u,v).
\end{equation*}
The weak formulation of problem \eqref{problem-neumann} is to find $u\in{\mathcal{H}}_J(\Omega)$ such that
\begin{equation*}
  \label{weak-neumann}
  {\mathcal{E}}(u,\varphi)=\int_\Omega f\varphi,
\end{equation*}
for every $\varphi\in{\mathcal{H}}_J(\Omega)$. Using the constant function $\varphi=1\in{\mathcal{H}}_J(\Omega)$ we get that a necessary condition to have a solution to problem \eqref{problem-neumann} is
\begin{equation}
  \label{compatibility}
  \int_\Omega f=0.
\end{equation}
Observe that for a constant it holds $\mathfrak{L}c=\mathcal{N}c=0$.
The following maximum principle is also immediate. Compare with Proposition~\ref{maximum}.
\begin{prop}\label{maximum-neumann}
  If $u\in{\mathcal{H}}_J(\Omega)$ satisfies
$\mathfrak{L}u\ge0$ in $\Omega$ and $\mathcal{N}u\ge0$ in $\Omega^c$
then $u$ is constant.
\end{prop}

Existence of solution is now proved using the compactness result obtained in Section~\ref{sect-inclusion}.

\begin{theorem}\label{th-Neumann-exist}
Assume  hypothesis \eqref{kernel2} with $\ell(0^+)=\infty$. Then given any $f\in L^2(\Omega)$ satisfying \eqref{compatibility} there exists a solution $u\in{\mathcal{H}}_J(\Omega)$ to problem~\eqref{problem-neumann}, unique up to additive constants.
\end{theorem}
\noindent{\it Proof.} Let $T_0\,:\,L^2(\Omega)\to{\mathcal{H}}_J(\Omega)$ be the operator defined by $T_0h=v$, where $v$ is the unique solution to the problem
$$
\begin{cases}
v+\mathfrak{L}v=h,& \mbox{in }\Omega,\\ \mathcal{N}v=0,& \mbox{in }\Omega^c.
\end{cases}
$$
The existence of such a solution follows from  Riesz representation Theorem.  Let $T\,:\,L^2(\Omega)\to L^2(\Omega)$ be defined by $Th=\left.T_0h\right|_\Omega$. Thanks to Theorem~\ref{th-compact2} this operator $T$ is compact, and it is also easily seen to be self-adjoint. Proposition~\ref{maximum-neumann} implies that $ker(I-T)$ consists only on constant functions. Therefore, for every $f\in (ker(I-T))^\perp$, that is, for every $f\in L^2(\Omega)$ with $\int_\Omega f=0$, there exists a function $w\in L^2(\Omega)$ satisfying $(I-T)w=f$. The function $u=T_0w$ satisfies
$$
\begin{cases}
u+\mathfrak{L}u=w,& \mbox{in }\Omega,\\ \mathcal{N}u=0,& \mbox{in }\Omega^c,
\end{cases}
$$
but in $\Omega$ it holds $w=f+Tw=f+T_0w=f+u$, so that $u$
solves problem~\eqref{problem-neumann}.
\qed

If the kernel does not decay too fast at infinity then any solution stabilizes  to a certain average. Let
\begin{equation*}\label{Jatinf}
  W(y)=\lim_{|x|\to\infty}\frac{J(x,y)}{J(x,0)}\qquad\text{for  } y\in\Omega,
\end{equation*}
and assume $0<c_1\le W(y)\le c_2<\infty$ for every $y\in\Omega$.

\begin{prop}
  \label{neumann-stabil}
  If $u$ is a solution to problem \eqref{problem-neumann}
then
\begin{equation*}
  \lim_{|x|\to\infty}u(x)=\frac{\displaystyle\int_\Omega u(y)W(y)\,dy}{\displaystyle\int_\Omega W(y)\,dy}.
\end{equation*}
\end{prop}
\noindent{\it Proof.} For every $\varepsilon>0$ there exists some $R>0$ such that
$$
(1-\varepsilon) J(x,0)W(y)<J(x,y)<(1+\varepsilon) J(x,0)W(y)
$$
for every $y\in\Omega$, $|x|>R$. Now the condition
$$
0=\mathcal{N}u(x)=
\int_{\Omega}(u(x)-u(y))J(x,y)\,dy
$$
outside $\Omega$ implies
$$
u(x)=\frac{\displaystyle\int_\Omega u(y)J(x,y)\,dy}{\displaystyle\int_\Omega J(x,y)\,dy},
$$
and therefore, for $|x|>R$,
$$
\frac{1-\varepsilon}{1+\varepsilon}\,\frac{\displaystyle \int_\Omega u(y)W(y)\,dy}{\displaystyle \int_\Omega W(y)\,dy}<u(x)<
\frac{1+\varepsilon}{1-\varepsilon}\,\frac{\displaystyle \int_\Omega u(y)W(y)\,dy}{\displaystyle \int_\Omega W(y)\,dy}.
$$
\qed

If $J$ decays at infinity like a power of $|x-y|$ then $W(y)=1$ and any solution stabilizes to its standard mean in $\Omega$,
$$
\lim_{|x|\to\infty}u(x)=\frac1{|\Omega|}\int_\Omega u(y)\,dy.
$$
This is not true if $J$ decays exponentially or even has compact support.

We finally  may  consider also the problem with nontrivial Neumann data
\begin{equation*}\label{problem-neumann2}
\begin{cases}
\mathfrak{L}u=f,& \mbox{in }\Omega,\\ \mathcal{N}u=g,& \mbox{in }\Omega^c,
\end{cases}
\end{equation*}
In that case we must assume that there exists some regular function $\psi$ such that $\mathcal{N}\psi=g$ in $\Omega^c$, something that is not clear. We then would obtain that the function $z=u-\psi$ satisfies the homogenous problem and we are reduced to the previous situation.

\section*{Acknowledgments}

Work supported by the Spanish project  MTM2014-53037-P.

\end{document}